\documentclass[11pt]{article}
\usepackage{amsfonts,amssymb,amsmath}

\usepackage[english]{babel}
\newcommand{\halmos}{\rule{1ex}{1.4ex}}

\makeatletter \@addtoreset{equation}{section} \makeatother

\newtheorem{ittheorem}{Theorem}
\newtheorem{itlemma}{Lemma}
\newtheorem{itproposition}{Proposition}
\newtheorem{itdefinition}{Definition}
\newtheorem{itremark}{Remark}
\newtheorem{itcorollary}{Corollary}

\newenvironment{theorem}{\addtocounter{equation}{1}
\begin{ittheorem}}{\end{ittheorem}}

\newenvironment{lemma}{\addtocounter{equation}{1}
\begin{itlemma}}{\end{itlemma}}

\newenvironment{proposition}{\addtocounter{equation}{1}
\begin{itproposition}}{\end{itproposition}}

\newenvironment{definition}{\addtocounter{equation}{1}
\begin{itdefinition}}{\end{itdefinition}}

\newenvironment{remark}{\addtocounter{equation}{1}
\begin{itremark}}{\end{itremark}}

\newenvironment{corollary}{\addtocounter{equation}{1}
\begin{itcorollary}}{\end{itcorollary}}

\newenvironment{proof}{\noindent {\em Proof}.\,\,\,}
{\hspace*{\fill}$\halmos$\medskip}

\newcommand{\beq}{\begin{eqnarray}}
\newcommand{\eeq}{\end{eqnarray}}

\newcommand{\be}{\begin{equation}}
\newcommand{\ee}{\end{equation}}

\newcommand{\bl}{\begin{lemma}}
\newcommand{\el}{\end{lemma}}

\newcommand{\br}{\begin{remark}}
\newcommand{\er}{\end{remark}}

\newcommand{\bt}{\begin{theorem}}
\newcommand{\et}{\end{theorem}}

\newcommand{\bd}{\begin{definition}}
\newcommand{\ed}{\end{definition}}

\newcommand{\bp}{\begin{proposition}}
\newcommand{\ep}{\end{proposition}}

\newcommand{\bc}{\begin{corollary}}
\newcommand{\ec}{\end{corollary}}

\newcommand{\bpr}{\begin{proof}}
\newcommand{\epr}{\end{proof}}

\newcommand{\bi}{\begin{itemize}}
\newcommand{\ei}{\end{itemize}}

\newcommand{\ben}{\begin{enumerate}}
\newcommand{\een}{\end{enumerate}}

\newcommand{\Z}{\mathbb Z}
\newcommand{\R}{\mathbb R}
\newcommand{\N}{\mathbb N}

\newcommand{\pee}{\mathbb P}

\newcommand{\gee}{\ensuremath{\mathcal{G}}}
\newcommand{\geeg}{\ensuremath{\mathcal{G}(\gamma)}}
\newcommand{\g}{\gamma)}
\newcommand{\geegu}{\ensuremath{\mathcal{G}(\gamma^\mu)}}

\newcommand{\s}{\ensuremath{\mathcal{S}}}
\newcommand{\loc}{\ensuremath{\mathcal{L}}}
\newcommand{\fe}{\ensuremath{\mathcal{F}}}

\newcommand{\la}{\ensuremath{\Lambda}}
\newcommand{\si}{\ensuremath{\sigma}}
\newcommand{\om}{\ensuremath{\omega}}
\newcommand{\Om}{\ensuremath{\Omega}}
\newcommand{\epsi}{\ensuremath{\epsilon}}
\newcommand{\m}{\ensuremath{\mathcal{M}_1^+}}
\newcommand{\minv}{\ensuremath{\mathcal{M}_{1,\rm{inv}}^+}}
\newcommand{\mee}{\ensuremath{\mathcal{M}}}

\def\a{\alpha}
\def\b{\beta}
\def\d{\delta}

\def\phi{\varphi}
\def\g{\gamma}
\def\l{\lambda}
\def\k{\kappa}

\def\r{\rho}

\def\x{\xi}

\def\D{\Delta}
\def\L{\Lambda}

\def\O{\Omega}

\def\del{\partial}

\def\GG{{\mathcal G}}
\def\HH{{\mathcal H}}

\def\bt{\tilde \b}

\def\1{1}

\def\N{{\Bbb N}}

\def\R{{\Bbb R}}

\def\Z{{\Bbb Z}}

\def\P{{\Bbb  P}}
\def\one{\hbox{J}\kern-.2em\hbox{I}}

\def\ba{{\backslash}}
\def\sb{{\subset}}

\def\emp{{\emptyset}}
\def\inn{\hbox{\rm int}}

\begin{document}

\title{{\bf Variational principle for generalized Gibbsian measures}}

\author{Christof K\"ulske\footnote{WIAS, Mohrenstrasse 39, 10117 Berlin, Germany.
  E-Mail: kuelske@wias-berlin.de}, Arnaud Le Ny\footnote{Eurandom, L.G. 1.48, TU Eindhoven, Postbus 513,
  5600 MB Eindhoven, The Netherlands. E-Mail:
  leny@eurandom.tue.nl}, Frank Redig\footnote{Faculteit Wiskunde En Informatica TU Eindhoven, Postbus 513,
  5600 MB Eindhoven, The Netherlands. E-Mail:f.h.j.redig@tue.nl}}

\maketitle



\maketitle

\footnotesize
\begin{center}
{\bf Abstract:} 
\end{center}

We study the thermodynamic formalism for generalized Gibbs measures, such as renormalization group transformations  of Gibbs measures or joint measures of disordered spin systems. We first show existence of the relative entropy density and obtain a familiar expression in terms of entropy and relative energy for "almost Gibbsian measures" (almost sure continuity of conditional probabilities). We also describe these measures as equilibrium states and establish an extension of the usual variational principle. As a corollary, we obtain a full variational principle for quasilocal measures. For the joint measures of the random field Ising model, we show that the weak Gibbs property holds, with an almost surely rapidly decaying translation invariant potential. For these measures we show that the variational principle fails as soon as the measures loses the almost Gibbs property. These examples suggest that the class of weakly Gibbsian measures is too broad from the perspective of a reasonable thermodynamic formalism.

\normalsize

\vspace{12pt}

\section{Introduction}
Since the discovery of the Griffiths-Pearce singularities
of renormalization group transformations \cite{GP,EFS}, it has been a challenging
question whether the classical Gibbs formalism can be
extended
in such a way as to incorporate renormalized low
temperature phases, so that renormalizing the measure can really be viewed as
a transformation on the level of Hamiltonians.
Later on, many other examples of ``non-Gibbsian" measures
appeared
in the context of joint measures of disordered spin
systems \cite{ku}, time
evolution of Gibbs measures \cite{EFHR}, and dynamical systems \cite{MRTV},
providing further
motivation for the construction of a generalized Gibbs
formalism.

As soon as the first examples of non-Gibbsian measures
appeared, Dobrushin proposed a program of ``Gibbsian restoration
of non-Gibbsian fields", arguing that the phenomenon of non-Gibbsianness 
is caused by ``exceptional" configurations, which are
negligible in the measure-theoretic sense. He thus proposed the
notion of a ``weakly Gibbsian" measure, where the existence of the
finite-volume Hamiltonian is not required uniformly in the
boundary condition, but only for boundary conditions in a set of
measure one. This is clearly enough to define the Gibbsian form of
the conditional probabilities, and Gibbs measures via the DLR equations. Since
Dobrushin and Shlosman (1997), many papers have been written showing the ``weak Gibbs"
property of renormalized low temperature phases, see e.g.,
\cite{BKL,MRV,MRSV,MVDV}, and of joint
measures of disordered spin systems, \cite{ku,ku1}.
Parallel to this, Fern\'andez and Pfister (1997) developed ideas
about generalized regularity properties of the conditional
probabilities. They proved that the decimation of the low
temperature plus phase of the Ising model is consistent with a
monotone right-continuous system of conditional probabilities. In
the framework of investigating regularity of the conditional
probabilities, the notion of ``almost Gibbs" was
introduced in \cite{MRV}. A measure $\mu$ is called almost Gibbs if its
conditional probabilities have a version which is continuous on a
set of $\mu$-measure one. If one does not insist on ``absolute"
convergence of the sums of potentials constituting finite-volume
Hamiltonians, then almost Gibbs implies weak Gibbs, but the
converse is not true, see \cite{lef2,MRV}. In \cite{FLNR}
it is proved that e.g. the decimation of the plus phase of the low
temperature Ising model is almost Gibbs, and the criterion to
characterize an essential point of discontinuity of the
conditional probabilities given in \cite{EFS} strongly suggests
that many other examples of renormalized low temperature phases
are almost Gibbs. The investigation of generalized Gibbs properties of the non-Gibbsian measures which appears e.g. as transformations of Gibbs measures is called the ``first part of the Dobrushin  program".

The  ``second part of the Dobrushin program" then consists in
building a thermodynamic formalism within the new class of
``generalized Gibbs measures". The question
whether in the context of weakly Gibbsian measures there
is a reasonable notion of ``physical equivalence" -- i.e., if two
systems of conditional probabilities share a Gibbs measure, then
they are equal --  is already raised in 
\cite{BKL}. In the classical Gibbs formalism, physical equivalence corresponds to zero relative entropy density, or zero ``information distance". Generally speaking, one would like to obtain a
relation between vanishing relative entropy density and
conditional probabilities. For Gibbs measures with a translation
invariant uniformly absolutely convergent potential, a translation
invariant probability measure $\mu$ has zero relative entropy
density $h(\mu|\nu)$ with respect to a Gibbs measure $\nu$ if and
only if $\mu$ is Gibbs with the same potential. Physically
speaking, this means that the only minimizers of the free energy
are the equilibrium phases. In complete generality, i.e. without
any locality requirements,
$h(\mu|\nu)=0$
does not imply that $\mu$ and $\nu$ have anything in common, see
e.g. the example in \cite{XI} where a  measure $\nu$ is constructed such that
for any translation invariant probability measure $h(\mu|\nu)=0$.

In this paper we investigate the relation between $h(\mu| \nu)=0$
and the property of having a common system of conditional
probabilities for general quasilocal measures, almost Gibbsian measures  and weakly Gibbsian measures. We will work in the context
of lattice spin systems with a single-site spin taking a finite
number of values. Let $\gamma$ denote a translation invariant
system of conditional probabilities, and $\gee_{\rm{inv}}
(\gamma)$ the set of all translation invariant probability
measures having $\gamma$ as a version of their conditional
probabilities. If $\gamma$ is continuous then, for
$\nu\in\gee_{\rm{inv}} (\gamma)$, we obtain $h(\mu| \nu)=0$ if and
only if $\mu\in\gee_{\rm{inv}} (\gamma)$. If $\gamma$ is
continuous $\mu$-almost everywhere, then we obtain that
$h(\mu|\nu)=0$ and $\nu\in\gee_{\rm{inv}} (\gamma)$ implies
$\mu\in\gee_{\rm{inv}} (\gamma)$. More generally, for $\nu\in\gee_{\rm{inv}}
(\gamma)$ and $\mu \in \m$  concentrating on a set of ``good configurations", we
obtain the existence of $h(\mu|\nu)$, an explicit
expression for it where $\nu$ enters only through its conditional
probabilities and the relation  $h(\mu|\nu)=0$ implies $\mu \in\gee_{\rm{inv}}
(\gamma)$. The ``good configurations" here are defined such
that a telescoping procedure - inspired by the method of Sullivan
\cite{Su} -  converges almost surely. These results, together with
some examples of non-Gibbsian measures
to which they apply show that almost Gibbsian
measures exhibit a reasonable thermodynamic formalism. The fact
that some concentration properties of the measures are required is
reminiscent of the situation in unbounded spin systems, an analogy already pointed out by
Dobrushin, \cite{Pirlot}.

The context of joint measures of disordered spin systems provides
a good source of examples for validity and failure of the
variational principle. Here by joint measure
we mean the joint distribution both of the spins and the disorder. In these examples (especially for the
random field Ising model) there is a precise criterion separating
the almost Gibbsian case from the weakly Gibbsian case. In
particular, for the random field Ising model, the joint measure is
always weakly Gibbs, and at low temperatures we prove here that it
even admits a translation invariant potential which decays almost
surely as a stretched exponential (so in particular converges {\em
absolutely} a.s.). If there is no phase transition, then the joint
measure for the random field Ising model is almost Gibbs (but
not Gibbs in dimension two at low temperature). In the almost Gibbsian
regime we obtain the variational principle whereas
in the weakly but not-almost Gibbsian regime we show the invalidity of
the variational principle. More precisely, in that case the joint measure for
the minus phase ($K^-$) is not consistent with the (weakly
Gibbsian) system of conditional probabilities of the plus phase
($K^+$), but one easily obtains that the relative entropy densities
$h(K^-|K^+)=h(K^+|K^-)=0$. Physically speaking, this means that we
are in the pathological situation where a minimizer of the free
energy is not a ``phase" (in the DLR sense). At the same time, we
also treat the joint measures in a very broad sense, i.e., for
possibly non-i.i.d. disorder, we prove existence of relative entropy
density, give an explicit representation in terms of the
defining potentials, and discuss implications of our results for the Morita approach \cite{Mo}.

Our paper is organized as follows: in Section 2 we
introduce basic definitions
and notations, discuss the different generalized Gibbs
measures and define
the variational principle. In Section 3 we prove the
variational
principle for some class of almost Gibbsian measures, using
the technique
of ``relative energies" of Sullivan \cite{Su}. In Section 4 we
prove the variational
principle for measures with translation invariant
continuous system of
conditional probabilities. In Section 5 we give the example
of the GriSing
random field and the decimation of the low-temperature plus
phase of the Ising model.
In Section 6 we discuss the examples of joint measures of
disordered spin systems.

\section{Preliminaries}

\subsection{Configuration space}

The configuration space is an infinite product space
$\Omega= E^{\Z^d}$ with $E$ a finite set.
Its Borel-$\si$-field is denoted by $\fe$.
 We denote by $\mathcal{S}=\big\{ \la \subset
\mathbb{Z}^d, |\la| < \infty \big\}$ the set of the finite subsets
of $\mathbb{Z}^d$ and for any $\la \in \s$,
$\Omega_\Lambda = E^\la$. $\fe_\la$ denotes the
$\sigma$-algebra generated by $\{ \si (x) :x\in\la \}$.
For all $\sigma,\omega \in
\Omega$, we denote $\sigma_{\la},\omega_{\la}$ the projections on
$\Omega_{\la}$ and also write $\sigma_{\la}\omega_{\la^c}$ for
 the configuration which agrees with $\sigma$ in $\Lambda$ and with
$\omega$ in $\la^c$. 
The set of  probability measures on $(\Omega,\fe)$
is denoted by $\m$. A function $f$ is said to be {\em local} if there exists $\Delta
 \in \s$ such that $f$ is $\mathcal{F}_{\Delta}$-measurable. We denote
by $\loc$ the set of all local functions. The uniform
closure of $\loc$ is $C(\Omega)$, the set of continuous functions
on $\Omega$.

On $\Omega$, translations
 $\{\tau_x:x\in\mathbb{Z}^d\}$
are defined via
$(\tau_x \omega)(y) =
\omega(x+y)$, and similarly on functions:
$\tau_x f(\omega)=
f(\tau_{x}\omega )$, and on measures
$\int f d\tau_{x}\mu=\int (\tau_{x}f) d\mu$.
The set of
translation invariant probability measures on $\Omega$
is denoted by $\minv$.

We also have a partial order: 
$\eta\leq\zeta$ if and only if for all $x\in\Z^d$, $\eta(x)\leq\zeta(x)$.
A function $f:\Omega\to\R$ is called monotone if
$\eta\leq\zeta$ implies $f(\eta)\leq f(\zeta)$. This order induces stochastic domination on $\m$: $\mu \; \preceq \; \nu$ if and only if $\mu(f) \leq \nu(f)$ for all $f$ monotone increasing.

\subsection{Specification and quasilocality}
\begin{definition}\label{chp:def1}
A {\em specification} on $(\Omega,\mathcal{F})$ is a family
$\gamma =\{\gamma_{\Lambda},\Lambda \in \mathcal{S}\}$ of
probability kernels from $\Omega_{\la^c}$ to $\mathcal{F}$ that
are
\begin{enumerate}
\item \emph{Proper}: For all $B \in
  \mathcal{F}_{\Lambda^{c}}$,
  $\gamma_{\Lambda}(B|\omega)=\mathbf{1}_{B}(\omega)$.
\item \emph{Consistent}: If $\Lambda \subset \Lambda'$ are finite
  sets, then $\gamma_{\Lambda'} \gamma_{\Lambda} = \gamma_{\Lambda'}$.
\end{enumerate}
\end{definition}
The notation $\gamma_{\Lambda'} \gamma_{\Lambda}$ refers to the
composition of probability kernels: for
$A\in\fe$, $\omega\in\Omega$:\[ (\gamma_{\Lambda'}
\gamma_{\Lambda})(A|\omega) =
\int_{\Omega}\gamma_{\Lambda}(A|\omega')\gamma_{\Lambda'}(d\omega'
| \omega). \] These kernels also acts on bounded measurable
functions $f$: \[
\gamma_{\la}f(\omega)=\int f(\si) \gamma_{\la}(d \si | \omega)\]
and on measures $\mu$: \[
\mu\gamma_{\la}(f)\equiv \int f d \mu \gamma_{\la}=\int ( \gamma_{\la}f ) d
\mu.\]
 A specification is a strengthening of the
notion of a system of proper regular conditional probabilities.
Indeed, in the former, the consistency condition (2) is required
to hold for \emph{every} configuration $\omega\in\Omega$, and not
only for almost every $\omega\in\Omega$.  This is because the
notion of specification is defined without any reference to a
particular measure.
A specification $\gamma$ is translation invariant if
for all $A\in\fe$, $\la\in\s$, $\omega\in\Omega$:
\[
\gamma_{\la +x} (A |\omega) = \gamma_\la (\tau_x A|\tau_x\omega)
\]
In this paper we will {\em always} restrict to the case
of non-null specifications, i.e., for any $\la\in\s$,
there exist $0< a_\la <b_\la <1$ such that
\[
a_\la < \inf_{\si,\eta} \gamma_\la (\si|\eta)\leq \sup_{\si,\eta} \gamma_\la
(\si|\eta)<b_\la.
\]

\bd

A probability measure $\mu$ on $(\Omega,\mathcal{F})$ is said to
be \emph{consistent} with a specification $\gamma$ (or specified
by $\gamma$) if the latter is a realization of its finite-volume
conditional probabilities, that is, if for all $A \in \mathcal{F}$
and $\la \in \s$, and for $\mu$-a.e. $\om$,

\be\label{chp:DLR}  \mu[A |\mathcal{F}_{\Lambda^{c}}](\,\om\,)=
\gamma_{\Lambda}(A |\,\om\,) .  \ee
\ed
Equivalently, $\mu$ is consistent with $\gamma$ if
\[ \int (\gamma_{\la}f ) d \mu =\int f d \mu
\]
for all $f\in C(\Omega)$.
We denote by $\geeg$ the set of measures
consistent with $\gamma$. For a translation
invariant specification, $\gee_{\rm{inv}}(\gamma)$ is
the set of
translation invariant elements of $\geeg$ .
\begin{definition}\label{chp:def2}
\begin{enumerate}
\item A specification $\gamma$ is \emph{quasilocal} if for each $\Lambda
\in \mathcal{S}$ and each $f$ local, $\gamma_{\Lambda}f\in C(\Omega)$.
 \item A probability measure $\mu$ is
\emph{quasilocal} if it is consistent with some quasilocal
specification.
\end{enumerate}
\end{definition}

\subsection{Potentials and Gibbs measures}
Examples of quasilocal measures are
{\em Gibbs measures}
defined via potentials.

\bd
\ben
\item
A {\em potential} is a family $\Phi=\{\Phi_A : A \in \s\}$ of local
functions such that for all $A \in \s$, $\Phi_A$ is
$\mathcal{F}_A$-measurable.
\item A potential is translation
invariant
if for all $A\in\s$, $x\in\Z^d$ and $\om\in\Omega$:
\[
\Phi_{A+x} (\om ) = \Phi_A(\tau_x\om)
\]
\een
\ed
\bd

A potential is said to be
\ben
\item
{\em Convergent
at} the configuration $\omega$ if
for all $\la\in\s$ the sum 
\be\label{chp:pot1}
\sum_{A \cap \la \neq \emptyset} \Phi_A(\omega)
\ee
is convergent.
\item {\em Uniformly convergent} if  convergence in
(\ref{chp:pot1}) is uniform in $\omega$.
\item
{\em Uniformly absolutely
convergent} (UAC) if for all $\la\in\s$
\[
\sum_{A \cap \la \neq \emptyset}
\sup_{\omega} |\Phi_A(\omega)| \; < \; \infty.
\]
\een
\ed
For a general potential $\Phi$, we define the measurable
set of its points of convergence:
\[
\Omega_\Phi = \{ \omega\in\Omega:\ \Phi \; \mbox{is convergent
at} \ \omega \}.
\]

In order to define Gibbs measures, we consider a UAC potential and define its
{\em  finite-volume Hamiltonian} for
$\la \in \s$ and boundary condition $\omega \in \Omega$
by
\[
H^{\Phi}_{\la}(\sigma | \omega)=
\sum_{A \cap \la \neq \emptyset}
\Phi_A(\sigma_{\la}\omega_{\la^c}).
\]

\bd\label{chp:def3} Let $\Phi$ be UAC. The
{\em Gibbs specification} $\gamma^\Phi$  with potential $\Phi$ 
is defined by
\begin{displaymath}
\gamma_{\la}^{ \Phi}(\si|\omega)=\frac{1}{Z^{ \Phi}_{\la}(\omega)}
e^{-H^{\Phi}_{\la}(\sigma|\omega)} 
\end{displaymath}
where the partition function $Z^{ \Phi}_{\la}(\omega)$ is the
normalizing constant. \ed A measure  $\mu$ is a  {\em Gibbs
measure} if there exists a UAC potential  $\Phi$  such that $\mu
\in \gee(\gamma^{\Phi})$. Gibbs measures are quasilocal and
conversely, any non-null quasilocal measure can be written in a
Gibbsian way (see \cite{K} and more details in Section 4). 

\subsection{Generalized Gibbs measures}

\begin{definition}\label{chp:def4}
A measure $\nu$ is {\em weakly Gibbs} if there exists a potential
$\Phi$ such that
$\nu (\Omega_\Phi )=1$ and
\begin{displaymath}
\nu\left[\sigma_\la|\fe_{\la^c}\right](\omega)= \frac{e^{-H^{
\Phi}_{\la}(\sigma|\omega)} }{Z^{ \Phi}_{\la}(\omega)}
\end{displaymath}
for $\nu$-almost every $\omega$.
\end{definition}

\br Some authors insist on the almost surely {\em absolute}
convergence of the sums defining $H^\Phi_\la$. 
However, for the definition of the weakly Gibbsian
specification
there is no reason to prefer absolute convergence.
\er
\begin{definition}\label{chp:def5}
Let $\gamma$ be a specification. A configuration $\omega$ is said to be a point of continuity for $\gamma$
if for all $\la \in \s$, $f \in \loc$, $\gamma_{\la}f$ is continuous at $\omega$.
\end{definition}
For a given $\gamma$, $\Omega_\gamma$ denotes
its measurable set of points of continuity.
\begin{definition}\label{chp:def6}
A measure $\nu$ is called {\em almost Gibbs} if there exists a
specification $\gamma$ such that $\nu \in \geeg$ and
$\nu(\Omega_{\gamma})=1$.
\end{definition}
If $\nu$ is almost Gibbs, then
there exists
an almost surely
convergent potential $\Phi$ such that $\nu$ is weakly Gibbsian for
$\Phi$, and thus almost Gibbsianness implies weak Gibbsianness. 
The converse is not true: a measure can be weakly Gibbs
and for the associated potential $\Phi$,
$\Omega_{\gamma^\Phi}$ is of measure zero \cite{lef2,MRV}.
If a measure is almost Gibbs and
translation invariant, then the corresponding potential
can be chosen to be translation invariant.

\subsection{Relative entropy and variational principle}
For
$\mu,\nu \in \minv$, the \emph{finite-volume relative
entropy} at volume $\Lambda \in \mathcal{S}$ of $\mu$ relative to
$\nu$ is defined as
\begin{equation}\label{eq:rex}
h_{\Lambda}(\mu | \nu)\;=\;\left\{
\begin{array}{ll}\; \displaystyle \int_{\Omega}
\frac{d\mu_{\Lambda}}{d\nu_{\Lambda}}
 \log\frac{d\mu_{\Lambda}}{d\nu_{\Lambda}}\, d \nu
 \; \; \; \textrm{if} \;   \;\mu_\Lambda \ll \nu_\Lambda \\[15pt]
+ \infty \; \; \; \; \; \; \; \; \; \; \; \;\; \; \; \; \; \; \;
\; \; \; \; \;\; \; \; \textrm{otherwise}.
\end{array} \right.
\end{equation}
The notation $\mu_\Lambda$ refers to the
distribution of $\omega_\la$ when $\omega$ is distributed
according to $\mu$. By Jensen's inequality, $h_\la (\mu|\nu) \geq 0$. 
The
\emph{relative entropy} of $\mu$ relative to $\nu$ is the limit
\begin{equation}\label{eq:red}
h(\mu | \nu)\;=\;\lim_{n \to \infty} \frac{1}{| \Lambda_n
|}h_{\Lambda_n}(\mu |
  \nu) 
\end{equation}
where $\Lambda_n = [\-n,n]^d\cap\Z^d$ is a sequence
of cubes (this can be replaced by a Van Hove sequence).
In what follows, if we write $\lim_{\la\uparrow\Z^d}f(\la)$
we
mean that the limit is taken
along a Van Hove sequence. The defining limit (\ref{eq:red}) is known to exist if $\nu\in \minv$
is a translation invariant Gibbs measure with a {\em translation
invariant} UAC potential
and $\mu\in
\mathcal{M}_{1,\rm{inv}}^{+}$ arbitrary.
The Kolmogorov-Sinai entropy
$h(\mu)$ is defined for $\mu\in\minv$:
\be\label{chp:KS}
h(\mu) = -\lim_{n\to\infty} \frac{1}{|\la_n|} \sum_{\si_{\la_n}}\mu(\si_{\la_n})
\log\mu\left(\si_{\la_n}\right).
\ee

We are now
ready to state the variational principle for specifications and
measures, which gives
a relation between zero relative entropy
and equality of conditional probabilities.

\begin{definition}\label{chp:def7}
 Let $\gamma$ be a
specification, $\nu \in \gee_{\rm{inv}}(\gamma)$ and $\mee\subset\minv$. We say that a
{\em variational principle holds for the triple $(\gamma,\nu,\mee)$}  if
\begin{description}
\item[(0)] $h(\mu|\nu)$ exists for all $\mu \in \mee$.
\item[(1)] $\mu \in \gee_{\rm{inv}}(\gamma)\cap\mee$ implies $h(\mu|\nu)=0$.
\item[(2)] $h(\mu|\nu)=0$ and $\mu \in \mee$ implies $\mu \in
\gee_{\rm{inv}}(\gamma)$.
\end{description}
\end{definition}
Items (1) and (2) are called the first and second part of the
variational principle. The second part is true for any translation
invariant quasilocal measure $\nu$ \cite{G} (with $\mee=\minv$). The first part is
proved for translation invariant Gibbs measures associated with a
translation invariant UAC potential (with $\mee=\minv$ also). We extend
this result to any translation invariant quasilocal measure in
Section 4. In \cite{FLNR}, the second  part has been proved for
some renormalized non-Gibbsian FKG measures.
In general, the set $\mee$ will be a set of translation invariant
probability measures concentrating on ``good configurations"
(e.g., points of continuity of conditional probabilities).

\section{Variational principle for generalized Gibbs measures}
We study the variational principle for generalized Gibbs measures.
We first prove the second
part for almost Gibbsian measures, which is a rather straightforward
technical extension of Georgii (1988), Chapter 15.
\subsection{Second part of the variational principle for almost
Gibbsian measures}
\begin{theorem}\label{chp:thm1}
Let $\gamma$ be a translation invariant specification on
$(\Omega,\mathcal{F})$ and $\nu \in \gee_{\rm{inv}}(\gamma)$. For
all $\mu \in \minv$,
\begin{displaymath}
\begin{array}{lll}
h(\mu | \nu)=0 \\
\mu(\Omega_{\gamma})=1
\end{array} \Bigl \} \Longrightarrow \; \mu \in \gee_{\rm{inv}}(\gamma)
\end{displaymath}
and thus such a measure $\mu$ is almost Gibbs w.r.t. $\gamma$.
\end{theorem}
\bpr Choose $\nu \in \gee_{\rm{inv}}(\gamma)$ and $\mu$ such that
$h(\mu | \nu)=0$. We have to prove that for any $g \in \loc,\la
\in \s$:
\be\label{blub}
\mu(\gamma_{\la}g -g)=0.
\ee
Fix $g \in \loc$ and $\Delta \in \s$ such that $g$ is $\fe_{\Delta}$-measurable.
 The hypothesis
\begin{equation}
h(\mu | \nu)=\lim_{\la \uparrow \Z^d} \frac{1}{|\la|} h(\mu | \nu)=0
\end{equation}
implies that for every $\la \in \s$, the density $f_{\la}=\frac{d
\mu_{\la}}{d \nu_{\la}}$ exists and is a bounded positive
$\fe_{\la}$-measurable function. Introduce local approximations
of $\gamma_{\la}g$:
\begin{eqnarray*}
g_n^-(\sigma)&=&\inf_{\omega \in \Omega} \gamma_{\la}g(\sigma_{\la_n}\omega_{\la_n^c})\\
g_n^+(\sigma)&=&\sup_{\omega \in \Omega}
\gamma_{\la}g(\sigma_{\la_n}\omega_{\la_n^c}).
\end{eqnarray*} 
In the quasilocal case, we have $g_n^+-g_n^- \rightarrow
0$ uniformly when $n$ goes to infinity, whereas here we have
$g_n^+-g_n^- \rightarrow 0$ on the set $\Omega_{\gamma}$ of
$\mu$-measure one, and hence,
by dominated convergence in $L^1 (\mu)$. To obtain
(\ref{blub}) decompose:
\begin{equation}
\mu(\gamma_{\la}g-g)=A_n +B_n+C_n+D_n
\end{equation}
where
\begin{eqnarray*}
A_n&=&\mu(\gamma_{\la}g-g_n^-)\\
B_n&=&\nu((g_n^- - \gamma_{\la}g) f_{\la_n \setminus \la})\\
C_n&=&\nu(f_{\la_n \setminus \la}(\gamma_{\la}g - g))\\
D_n&=&\nu((f_{\la_n \setminus \la}-f_{\la_n})g).
\end{eqnarray*}
Using
\begin{displaymath}
0 \leq \gamma_{\la}g - g_n^- \leq g_n^+ - g_n^-
\end{displaymath}
$A_n \rightarrow 0$ as $n$ goes to infinity.
For $B_n$, use
\begin{displaymath}
0 \leq |B_n| = \nu\left((\gamma_{\la}g - g_n^-)f_{\la_n \setminus \la}\right) \leq \nu(f_{\la_n \setminus \la}(g_n^+-g_n^-))=\mu(g_n^+-g_n^-),
\end{displaymath}
to obtain $B_n\to 0$ as $n\to\infty$.

Since $\nu\in\geeg$, and $f_{\la_n\setminus\la}\in\fe_{\la^c}$,
$C_n=0$.
The fact that $D_n\to 0$ follows from the assumption
of zero relative entropy density: see Georgii (1988), p 324.

\epr
\br
\ben
\item The role of $\mee$ in Definition \ref{chp:def7}
is played here by the set of measures concentrating
on the points of continuity  of $\gamma$ ($\mu\in\mee$ if
and only if $\mu (\Omega_\gamma ) =1$).
\item Remark that in Theorem 3.1, we do not ask any
concentration properties of $\nu$.
\een
\er
\subsection{First part of the variational principle for some almost
Gibbsian measures}

To obtain the first part of the variational
principle, it will turn out
that concentration of $\mu$ on the set $\Omega_\gamma$
is not the right condition. We need that some particular
class of ``telescoping configurations" are points
of continuity of the specification. This reminds of
asking continuity properties of the one-sided
conditional probabilities. In the case of
(uniformly) continuous specifications, this distinction
between one-sided and two-sided is of course not visible.

We choose a particular value written $+1$ in the state
space $E$ and denote by "+" the configuration whose value is $+1$
everywhere.  To any configuration  $\sigma \in \Omega$, we
associate the configuration $\si^+$  defined by
\begin{displaymath}
\sigma^+(x)=\left\{
\begin{array}{ll}
\sigma(x) \; \; \rm{if} \; \; \;x \; \leq \; 0\\
\\
+ 1\; \; \; \; \rm{if} \; \; \; x \; >  \; 0.
\end{array}
\right.
\end{displaymath}
Here, the order $\leq$ is lexicographic.
We define then $\Omega_{\gamma}^{<0}$ to be the subset of $\Omega$
of the configurations $\sigma$ such that the new configuration
$\si^+$ is a good configuration for $\gamma$:
\begin{displaymath}
\Omega_{\gamma}^{<0}=\bigl \{\sigma \in \Omega,\; \sigma^+
\in \Omega_{\gamma} \bigr \}.
\end{displaymath}
This set will be described in different examples in Section 5.

\subsubsection{Results} We consider a pair $(\gamma,\nu)$ with  $\nu \in \gee_{\rm{inv}}(\g)$ and a measure $\mu$ which
satisfies the
following condition:\\

{\bf Condition C1} \[ 
\mu(\Omega_{\gamma}^{<0})=1. \] \medskip
We also introduce 
\[
e_{\nu}^+:=-\lim_{\la \uparrow \Z^d}
\frac{1}{|\la|}
 \log \nu(+_{\la})
\]
whenever it exists.
\begin{theorem}\label{chp:thm2}
Under the  condition C1:
\begin{enumerate}
\item $h(\mu|\nu)$  exists if and only if $e_{\nu}^+$ exists and then
\begin{equation}\label{chp:exist}
h(\mu|\nu)= e_{\nu}^+ - h(\mu) - \int_{\Omega} \log
\frac{\gamma_0(\sigma^+ | \sigma^+)}{\gamma_0( + |
\sigma^+)}\mu(d\sigma).
\end{equation}
where $h(\mu)$ is the Kolmogorov-Sinai entropy of $\mu$.

\item If moreover $\mu \in \gee_{\rm{inv}}(\gamma)$ and
$e_{\nu}^+$ exists, then
\be\label{tension}
h(\mu |\nu)=\lim_{\la \uparrow \Z^d} \frac{1}{|\la|} \log
\frac{\mu(+_{\la})}{\nu(+_{\la})}.
\ee
\end{enumerate}
\end{theorem}
To get the more usual expression of the variational principle, we
add an extra condition to the condition C1: \\

{\bf Condition C2}
\[
\mu \in \mathcal{G}_{\rm{inv}}(\gamma) \; \rm{is \; such \; that
\;} \lim_{\la \uparrow \Z^d} \frac{1}{|\la|} \log
\frac{\mu(+_{\la})}{\nu(+_{\la})}=0.
\]

\begin{theorem}\label{chp:thm3}
Assume that conditions C1 and C2 are fulfilled. Then
\begin{enumerate}
\item $h(\mu |
\nu)=0$.
\item  $e_{\nu}^+$ exists and
$e_{\nu}^+=e_{\mu}^+$.
\item $h(\alpha|\nu)$ exists for all $\alpha \in \minv$ satisfying C1.
\end{enumerate}
\end{theorem}
As a corollary of these theorems, we obtain the usual first part
of the variational principle.
\begin{theorem}\label{chp:thm4}
Let $\mu \in \minv$ and $\nu \in \gee_{\rm{inv}}(\gamma)$ such
that conditions C1 and  C2 hold and $e_{\nu}^+$ exists. Then
\begin{enumerate}
\item $h(\mu|\nu)$ exists.
\item $\mu \in \gee_{\rm{inv}}(\gamma)$ implies $h(\mu|\nu)=0$.
\end{enumerate}
\end{theorem}

\br The existence of the limit defining $e_{\nu}^+$ is guaranteed for
e.g. renormalization
group transformations of Gibbs measures, and for $\nu$ with
positive correlations (by subadditivity). Moreover, in the case of
transformations of Gibbs measures, condition C2 is also easy to
verify. See Section 5 below. \er
\subsection{Proofs}
\begin{description}
\item[Proof of Theorem \ref{chp:thm2}]
\end{description}
First we need the following
\begin{lemma}\label{chp:lem1}
If $\mu(\Omega_{\gamma}^{<0})=1$, then
\ben
\item Uniformly in $\om \in \Omega$,
\begin{displaymath}
\lim_{n \to \infty} \frac{1}{|\la_n|} \int_{\Omega} \log
\frac{\gamma_{\la_n}(\si|\omega)}{\gamma_{\la_n}(+|\omega)} \mu(d
\sigma)=\int_{\Omega}\log
\frac{\gamma_0(\si^+|\si^+)}{\gamma_0(+|\si^+)} \mu(d\sigma).
\end{displaymath}
\item For $\nu\in\geeg$,
\begin{displaymath}
\lim_{n \to \infty} \frac{1}{|\la_n|} \int_{\Omega} \log
\frac{\nu(\si_{\la_n})}{\nu(+_{\la_n})} \mu(d \sigma)=
\int_{\Omega} \log \frac{\gamma_0(\si^+|\si^+)}{\gamma_0(+|\si^+)}
\mu(d\sigma).
\end{displaymath}
In particular, the limit depends only on the pair $(\gamma,\mu)$.
\een
\end{lemma}

\br

If $\mu$ is ergodic under translations, we have a
slightly stronger statement for item 1 : $\frac{1}{|\la_n|}
\int_{\Omega} \log
\frac{\gamma_{\la}(\si|\omega)}{\gamma_{\la}(+|\omega)} \mu(d
\sigma)$ converges in $\mathbb{L}^1(\mu)$ to $\int_{\Omega} \log
\frac{\gamma_0(\si^+|\si^+)}{\gamma_0(+|\si^+)} \mu(d\sigma)$,
uniformly in $\omega \in \Omega$.
\er
 \bpr
\begin{enumerate}
\item The proof uses  relative energies as in Sullivan (1973). For all $\la \in
\s,\sigma, \om \in \Omega$, we define,
\begin{center}
$E_{\la}^+(\sigma | \om)=\log
\frac{\gamma_{\la}(\sigma|\om)}{\gamma_{\la}(+|\om)} \; \;$
and $\; \; D(\si)= E^+_{\{ 0 \}} (\si|\si)=\log
\frac{\gamma_0(\si|\si)}{\gamma_0(+|\si)}$.
\end{center}
We consider an approximation  of $\sigma^+$ at finite volume $\la$ with boundary
condition $\om$ and define the {\em telescoping
configuration} $T_{\la}^{\om}[x,\si,+]$:
\[
T_{\la}^{\om}[x,\si,+](y)= \left\{
\begin{array}{lll}
\omega(y) \;  &\rm{if}& \; y \in \la^c\\
\\
\sigma(y) \;&\rm{if}& \; y \; \leq \; x, \; y \in \la\\
\\
+ 1 \; &\rm{if}& \;  y\; > \; x,\; y \in \la.\\
\\

\end{array}
\right.
\]
Using the consistency property of $\gamma$, we have by
telescoping, 
\begin{displaymath}
E^+_{\la}(\sigma|\om)=\sum_{x \in \la} E_x^+(\sigma |
T_{\la}^{\om}[x,\si,+]).
\end{displaymath}
By translation invariance of $\gamma$,
\begin{displaymath}
E_{\la}^+(\sigma | \om)=\sum_{x \in \la}
D(\tau_{-x}T_{\la}^{\om}[x,\si,+]).
\end{displaymath}

By translation invariance of $\mu$,
\begin{eqnarray*}
\int_{\Omega} E^+_{\la_n}(\sigma|\omega) \mu(d \sigma)&=&\sum_{x \in
\la_n} \int_{\Omega} D(\tau_{-x} T_{\la}^{\om}[x,\tau_x \si,+]) \mu(d \sigma).
\end{eqnarray*}
Therefore, we have to prove that, uniformly in $\om$,
\[
\lim_{n \to \infty} \frac{1}{|\la_n |} \left(\sum_{x \in
\la_n}\int_{\Omega}\Bigl[D(\tau_{-x} T_{\la_n}^{\om}[x,\tau_x
\si,+])-D(\si^+) \Bigr] \mu(d\sigma) \right) =0.
\]
By definition,
\begin{displaymath}
\tau_{-x}T_{\la_n}^{\om}[x,\tau_x \si,+]= \left\{
\begin{array}{lll}
\tau_{-x} \om(y) \;  &\rm{if}& \; y+x \in \la_n^c\\
\\
+ \; &\rm{if}& \; 0 \; < y\;, y+x \in \la_n\\
\\
 \sigma(y) \;&\rm{if}& \; y \leq \; 0, \; y+x \in \la_n.
\end{array}
\right.
\end{displaymath}

Now, pick $\epsilon > 0$, $\om \in \Omega$ and  $\sigma \in
\Omega^{< 0}_{\gamma}$. Using the fact that $\sigma^+$ is a point of
continuity of $D$, we choose $n_0$ such that $\xi|_{\la_{n_0}}=\si^+|_{\la_{n_0}}$ implies
$|D(\xi)-D(\si^+)|\leq \epsilon$. Remark
that $\tau_{-x}T_{\la_n}^{\om}[x,\tau_x \si,+]$ and $\si^+$ 
differ only on the set $\{ y\in\Z^d: x+y\in\la_n^c\}$. Therefore,
the difference $|D(\si^+)-D(\tau_{-x}T_{\la_n}^{\om}[x,\tau_x \si,+])|$
can only be bigger than $\epsilon$ for $x$ such that
$(\la_{n_0} -
x) \cap \Lambda_n^c \neq \emptyset$.

Therefore,
\begin{eqnarray*}
 \frac{1 }{| \la_n  |} \Bigl|\; \sum_{x \in \la_n} \big[ D(\tau_{-x}
T_{\la_n}^{\om}[x,\tau_x \si,+])-D(\sigma^+ \big]\; \Bigr|\\\leq
\epsilon + \; \frac{2}{| \la_n  |} \mid \mid D \mid \mid_{\infty}\; \Big| \{x \in \la_n\; : \; (\la_{n_0} -
x) \cap \Lambda_n^c \neq \emptyset \; \} \Big| 
\end{eqnarray*}
and this is less than $2 \epsilon$ for $n$ big
enough. So we obtain that
\begin{displaymath}
\frac{1}{| \la_n  |} \Bigl|\; \sum_{x \in \la_n} \big[ D(\tau_{-x}
T_{\la_n}^{\om}[x,\tau_x \si,+])-D(\sigma^+) \big] \; \Bigr|
\end{displaymath}
converges to zero on the set of $\Omega^{< 0}_{\gamma}$ of full
$\mu$-measure, uniformly in $\om$. By dominated convergence, we
then obtain
\begin{displaymath}
\lim_{n \to \infty} \sup_{\om} \frac{1}{| \la_n|}\int_{\Omega}
\Big |\sum_{x \in \la_n} \big[ D(\tau_{-x} T_{\la}^{\om}[x,\tau_x
\si,+])-D(\sigma^+) \bigr] \Big |\mu(d\sigma)=0
 \end{displaymath}
which implies statement 1 of the lemma.
\item Denote
\[
F_{\la_n}(\mu,\nu)=\frac{1}{ | \la_n |} \int_{\Omega} \log
\frac{\nu(\sigma_{\la_n})}{\nu(+_{\la_n})} \mu(d \sigma).
\]
Using $\nu \in \geeg$, we obtain
\begin{displaymath}
 F_{\la_n}(\mu,\nu)=\frac{1}{| \la_n  |}
\int_{\Omega} \log \frac{\int_{\Omega} \gamma_{\la_n}(\sigma |
\om) \nu(d \om)}{\int_{\Omega} \gamma_{\la_n}( + | \om) \nu(d
\om)} \mu(d \si).
\end{displaymath}
Use
\begin{displaymath}
\inf_{\om \in \Omega}\frac{\gamma_{\la_n}(\sigma | \om)}{
\gamma_{\la_n}( + | \om) } \; \leq \; \frac{\int_{\Omega}
\gamma_{\la_n}(\sigma | \om) \nu(d \om)}{\int_{\Omega}
\gamma_{\la_n}( + | \om) \nu(d \om)} \; \leq \; \sup_{\om \in
\Omega} \frac{\gamma_{\la_n}(\sigma | \om)}{ \gamma_{\la_n}( + |
\om)}.
\end{displaymath}
Let $\epsi >0$ be given and $\om=\om(n,\sigma,
\epsilon)$, $\om'=\om'(n,\sigma, \epsilon)$ such that
\[
\int_{\Omega} \inf_{\om \in \Omega}\log
\frac{\gamma_{\la_n}(\sigma | \om)}{ \gamma_{\la_n}( + | \om) }
\mu(d \sigma) \geq \int_{\Omega} \log
\frac{\gamma_{\la_n}(\sigma|\om(n,\sigma,
\epsilon))}{\gamma_{\la_n}(+|\om(n,\sigma, \epsilon))} - \epsi
\]
and
\[
\int_{\Omega} \sup_{\om \in \Omega}\log
\frac{\gamma_{\la_n}(\sigma | \om)}{ \gamma_{\la_n}( + | \om) }
\mu(d \sigma) \leq \int_{\Omega} \log
\frac{\gamma_{\la_n}(\sigma|\om'(n,\sigma,
\epsilon))}{\gamma_{\la_n}(+|\om'(n,\sigma, \epsilon))} + \epsi.
\]
Now use the first item of the lemma and choose $N$ such that for
all $n \geq N$,

\[
\sup_{\om} \Big| \frac{1}{|\la_n |}\int_{\Omega} \log \frac{\gamma_{\la_n}(\sigma |
\om) }{\gamma_{\la_n}( + | \om) } \mu(d \sigma) - \int_{\Omega}
D(\sigma^+) \mu(d \si) \Big| \leq \epsilon.
\]
For $n \geq N$, we obtain
\begin{displaymath}
\int_{\Omega} D(\sigma^+) \mu(d \sigma) -2 \epsi \; \leq \;
F_{\la_n} (\mu | \nu) \; \leq \; \int_{\Omega} D(\sigma^+)
\mu(d\si) + 2 \epsi.
\end{displaymath}
\end{enumerate}

\epr

{\bf Proof of Theorem \ref{chp:thm2}}

\begin{enumerate}
\item Denote
\begin{displaymath}
h_n(\mu | \nu):=\frac{1}{|\la_n|}\sum_{\sigma_{\la_n}}
\mu(\sigma_{\la_n})\log
\frac{\mu(\sigma_{\la_n})}{\nu(\sigma_{\la_n})}.
\end{displaymath}
We recall that for $\mu \in \mathcal{M}_{1,\rm{inv}}^+(\Omega)$,
the limit of $h_n(\mu):=-\frac{1}{|\la_n|}\sum_{\sigma_{\la_n}}
\mu(\sigma_{\la_n})\log \mu(\sigma_{\la_n})$ is the {\em Kolmogorov-Sinai entropy}
of $\mu$ denoted $h(\mu)$. We write
\[
h_n(\mu | \nu)=\, -h_n(\mu) -\,
\frac{1}{|\la_n|}\sum_{\sigma_{\la_n}} \mu(\sigma_{\la_n})\log
\frac{\nu(\sigma_{\la_n})}{\nu(+_{\la_n})}\, -
\,\frac{1}{|\la_n|}\log \nu(+_{\la_n}).
\]
When condition C1 holds, the asymptotic behavior of the second
term of the r.h.s. is given by  Lemma \ref{chp:lem1}. Hence, 
  the relative entropy exists if and only if $e_{\nu}^+$ exists,
  and it is given by (\ref{chp:exist}).
\item We consider $\mu \in \mathcal{G}_{\rm{inv}}(\gamma)$ such that
$\mu(\Omega_{\gamma}^{< 0})=1$ and use the following decomposition
of the finite-volume relative entropy:

\beq\label{chp:efv2}
& &h_n(\mu|\nu)= \\ 
& &\frac{1}{|\la_n|}\sum_{\sigma_{\la_n}} \mu(\sigma_{\la_n})\log
\frac{\mu(\sigma_{\la_n})}{\mu(+_{\la_n})} -
\frac{1}{|\la_n|}\sum_{\sigma_{\la_n}} \mu(\sigma_{\la_n})\log
\frac{\nu(\sigma_{\la_n})}{\nu(+_{\la_n})} + \frac{1}{|\la_n|}\log
\frac{\mu(+_{\la_n})}{\nu(+_{\la_n})}.\nonumber
\eeq
By Lemma \ref{chp:lem1}, in the 
limit $n\to\infty$, 
the first two terms of the
r.h.s. are functions of $\gamma$ rather than functions of  $\mu, \nu
\in \gee_{\rm{inv}}(\gamma)$ and cancel out. Hence, the relative
entropy exists if and only if the third term converges. Using Item
1 (existence of relative entropy), we obtain the existence of the
limit (\ref{tension}) and the equality
\[
h(\mu | \nu)=\lim_{n \to \infty} \frac{1}{| \, \la_n \, |}
\log \frac{\mu(+_{\la_n})}{\nu(+_{\la_n})}.
\]
\end{enumerate}

\begin{description}
\item[Proof of Theorem \ref{chp:thm3}]
\end{description}
\begin{enumerate}
\item This is direct consequence of Theorem \ref{chp:thm2} and (\ref{tension}): under the
conditions $C1$ and $C2$,  $h(\mu|\nu)=0$.

\item  The existence
of the relative entropy proves that $e_{\nu}^+$ exists and is
given by
\begin{displaymath}
e_{\nu}^+=h(\mu) + \int \log \frac{\gamma_0(\sigma^+
|\sigma^+)}{\gamma_0(+ | \si^+)} \mu(d \si).
\end{displaymath}
Combined with C2 this proves $e_{\mu}^+=e_{\nu}^+$.

\item  Consider any other
measure $\alpha \in \minv$ such that C1 holds. The existence of
the relative entropy $h(\alpha|\mu)$ follows by combining 
the existence of $e_{\nu}^+$ with Theorem \ref{chp:thm2}, and
\begin{displaymath}
h(\alpha|\nu)=e_{\nu}^+ - h(\alpha) - \int \log
\frac{\gamma_0(\sigma^+ |\sigma^+)}{\gamma_0(+ | \si^+)} \alpha(d
\si).
\end{displaymath}
If moreover $\alpha$ satisfies C2, we also obtain that
$e_{\alpha}^+$ exists and equals $e_{\nu}^+$.
\end{enumerate}

\subsection{Generalization}

In the hypothesis of the theorems above, the plus-configuration
plays a particular role of telescoping reference configuration.
Without too much effort, we obtain the following generalization
where we telescope w.r.t a random configuration $\xi$ chosen from
some translation invariant measure $\lambda$. Results of the
previous section are recovered by choosing $\lambda = \delta_+$.
The generalization to a random telescoping configuration will be
natural in the context of joint measures of disordered spin systems in Section 6.

For any $\xi,\si \in \Om$, we define the concatenated
configuration $\si^{\xi}$:
\[
\forall x \in \Z^d, \sigma^{\xi}(x)=\left\{
\begin{array}{ll}
\sigma(x) \; \; \rm{if} \; \; \;x \; \leq \; 0\\
\\
\xi(x)\; \;  \rm{if} \; \; \; x \; > \; 0.

\end{array}
\right.
\]
and the set $\Omega_{\gamma}^{\xi,<0}$ to be the subset of $\Omega
\times \Om$ of the configurations $(\sigma, \xi)$ such that the
new configuration $\si^{\xi}$ is a good configuration for
$\gamma$:
\begin{displaymath}
\Omega_{\gamma}^{\xi,<0}=\bigl \{(\sigma,\xi) \in \Omega \times
\Om,\; \sigma^{\xi} \in \Omega_{\gamma} \bigr \}.
\end{displaymath}
We also generalize $e_{\nu}^+$ and denotes
\be \label{elmu}
e_{\nu}^{\lambda}=-\lim_{\la \uparrow \Z^d }\frac{1}{|\la|} \int_{\Om} \log
\nu(\xi_{\la}) \lambda(d\xi)
\ee
 provided this limit
exists.

We consider a specification $\gamma$, measures $\nu \in
\gee_{\rm{inv}}(\gamma)$, $\mu,\lambda \in \minv$, and the following
conditions:
\begin{description}
\item[C'1] $\lambda \otimes \mu (\Omega_{\gamma}^{\xi,<0})=1$.
\item[C'2] $\lim_{\la \uparrow \Z^d}
\frac{1}{|\la|} \int_{\Om} \big( \log \frac{d \mu_{\la}}{d
\nu_{\la}} \big) (\xi_{\la})  \lambda(d \xi_{\la})=0$.
\end{description}
The following theorems are the straightforward generalizations of
Theorem \ref{chp:thm2} and \ref{chp:thm4}, and their proofs follow
the same lines.

\begin{theorem}\label{chp:thm5}
Under the  condition C'1,
\begin{enumerate}
\item $h(\mu|\nu)$  exists if and only if $e_{\nu}^{\lambda}$ exists and then
\begin{equation}\label{chp:rella}
h(\mu|\nu)=e_{\nu}^{\lambda}-h(\mu)  - \int_{\Om \times \Om} \log
\frac{\gamma_0(\sigma^{\xi}|\sigma^{\xi})}{\gamma_0(\xi |
\sigma^\xi)} \mu(d \si) \lambda(d \xi).
\end{equation}

\item If moreover $\mu \in \gee_{\rm{inv}}(\gamma)$ and
$e_{\nu}^{\lambda}$ exists, then
\begin{displaymath}
h(\mu | \nu)=\lim_{\la \uparrow \Z^d} \frac{1}{|\la|} \int_{\Om}
\big( \log \frac{d \mu_{\la}}{d \nu_{\la}} \big)(\xi_{\la})
\lambda(d \xi_{\la}).
\end{displaymath}
\end{enumerate}
\end{theorem}

\begin{theorem}\label{chp:thm6} Consider $\mu \in \minv, \gamma$ a specification, $\nu \in
\mathcal{G}_{\rm{inv}}(\gamma)$ such that $e_{\nu}^{\lambda}$
exists and conditions C'1 and C'2 are true. Then
\begin{enumerate}
\item  $h(\mu|\nu)$
exists and is given by (\ref{chp:rella}).

\item
$\mu \in \gee_{\rm{inv}}(\gamma)$ implies
$ h(\mu |\nu)=0$.
\end{enumerate}
\end{theorem}

\section{Variational principle for quasilocal measures}

The usual way to prove $\mu \in \mathcal{G}_{\rm{inv}}(\gamma)
\Longleftrightarrow h(\mu | \nu)=0$ in the Gibbsian context uses
that
 $\gamma$ is a specification associated with a translation invariant and UAC potential $\Phi$,
 and goes via existence and boundary condition independence of pressure (see \cite{G}).
Since for a general quasilocal specification $\gamma$ we cannot rely on the existence of such a potential (see \cite{K} and open problem in \cite{EFS}),
we show here that the weaker property of uniform convergence
of the vacuum potential which
can be associated to the quasilocal $\gamma$
(see \cite{K}) suffices to obtain zero relative entropy.

\begin{theorem}\label{chp:thm7}
Let $\gamma$ be a translation invariant quasilocal specification,
$\nu \in \gee_{\rm{inv}}(\gamma)$ and $\mu \in \minv$. Then $h(\mu
| \nu)$ exists for all $\mu \in \minv$ and
\[
\mu \in \mathcal{G}_{\rm{inv}}(\gamma) \, \Longleftrightarrow \,
h(\mu|\nu)=0.
\]
\end{theorem}
\bpr The implication of the left ({\em the  second part}) is
proved in \cite{G}. To prove the first part, we need the following
lemma to check hypothesis of Theorem \ref{chp:thm4}. Condition
C2 is trivially true when $\gamma$ is quasilocal
($\Omega_{\gamma}^{<0}=\Omega$).

\begin{lemma}\label{chp:lem2}
For all $\mu,\nu \in \gee_{\rm{inv}}(\gamma)$ with $\gamma$
translation invariant and quasilocal, $e^+_{\nu},e^+_{\mu}$ exist
and
\begin{displaymath}
\lim_{n \to \infty} \frac{1}{|\la_n|} \log
\frac{\mu(+_{\la_n})}{\nu(+_{\la_n})}=0.
\end{displaymath}
\end{lemma}

\bpr
Kozlov (1974) proves that to any translation
invariant quasilocal specification $\gamma$ there corresponds a translation
invariant uniformly convergent vacuum potential $\Phi$ such that
$\gamma=\gamma^\Phi$.

By uniform convergence, we have
\begin{equation}\label{chp:us}
\lim_{\la \uparrow \mathbb{Z}^d} \sup_{\sigma} \; \Bigl| \sum_{A \ni
0, A \cap \la^c \neq \emptyset} \Phi_A(\sigma) \Bigr| \; = 0.
\end{equation}
Remark that in (\ref{chp:us}), the absolute value is \emph{outside} the sum,
i.e., (\ref{chp:us}) means that the series $\sum_{A \ni 0}
\Phi_A(\sigma)$ is convergent in the sup-norm topology on
$C(\Omega)$, but not necessarily {\em absolutely
convergent}. We can define a
Hamiltonian and a partition function for any $\la \in
\mathcal{S},\eta,\sigma \in \Omega$, as usual:
\begin{equation}\label{chp:eqvac}
H_{\la}^{\eta}(\sigma)=\sum_{A \cap \la \neq \emptyset}
\Phi_A(\sigma_{\la}\eta_{\la^c})\; \; \; \; {\rm and} \; \; \; \;
Z_{\la}(\omega)=\sum_{\sigma \in \Omega}
e^{-H_{\la}^{\omega}(\sigma)}.
\end{equation}
Lemma \ref{chp:lem2} is now a direct consequence of the following

\begin{lemma}\label{chp:lem3}

\begin{enumerate}
\item
\begin{equation}\label{chp:e4.3.1}
 \lim_{n \to \infty} \sup_{\omega,\eta,\sigma}
\frac{1}{|\la_n|} \bigl|  H_{\la_n}^{\eta}(\sigma) -
H_{\la_n}^{\omega}(\sigma) \bigr| =0.
\end{equation}
\item
\begin{equation}\label{chp:e4.3.2}
\lim_{n \to \infty} \sup_{\omega,\eta}\frac{1}{|\la_n|} \log
\frac{Z_{\la_n}(\omega)}{Z_{\la_n}(\eta)} =0.
\end{equation}
\end{enumerate}
\end{lemma}
\begin{proof}
We follow the standard line of the argument used by Israel (1986) to
prove existence and boundary condition independence of the
pressure for a UAC potential, but we detail it because the vacuum
potential is only uniformly convergent. Clearly, (\ref{chp:e4.3.1})
implies (\ref{chp:e4.3.2}): for all $n \in \mathbb{N}$,
\begin{displaymath}
\exp\Big\{-\sup_{\omega,\eta,\sigma} \Bigl|
H_{\la_n}^{\eta}(\sigma) - H_{\la_n}^{\omega}(\sigma) \Bigr|
\Big\} \leq
\sup_{\omega,\eta}\frac{Z_{\la_n}(\omega)}{Z_{\la_n}(\eta)} \leq
\exp\Big\{\sup_{\omega,\eta,\sigma} \Bigl|
H_{\la_n}^{\eta}(\sigma) - H_{\la_n}^{\omega}(\sigma) \Bigr|
\Big\}.
\end{displaymath}
To prove (\ref{chp:e4.3.1}), we write
\begin{displaymath}
H_{\la_n}^{\eta}(\sigma) - H_{\la_n}^{\omega}(\sigma)=\sum_{A \cap
\la_n \neq \emptyset, A \cap \la_n^c \neq \emptyset} \Bigl[
\Phi_A(\sigma_{\la_n}\eta_{\la_n^c})-\Phi_A(\sigma_{\la_n}\omega_{\la_n^c})
\Bigr].
\end{displaymath}
and we first remark:
\begin{displaymath}
\frac{1}{|\la_n|} \Bigl| \sum_{A \cap \la_n \neq \emptyset, A \cap
\la_n^c \neq \emptyset} \Bigl[
\Phi_A(\sigma_{\la_n}\eta_{\la_n^c})-\Phi_A(\sigma_{\la_n}\omega_{\la_n^c})
\Bigr] \Bigr| \leq \frac{2}{|\la_n|} \sum_{x \in \la_n} \sup_{\sigma}
\Bigl| \sum_{A \ni x, A \cap \la_n^c \neq \emptyset}
\Phi_A(\sigma)\Bigr|.
\end{displaymath}
We obtain
\begin{eqnarray*}
\sup_{\sigma} \Bigl| \sum_{A \ni x, A \cap \la_n^c \neq \emptyset}
\Phi_A(\sigma)\Bigr| &=& \Bigl| \sum_{A \ni x} \Phi_A(\sigma)
\;-\sum_{A \ni x, A \subset
\la_n} \Phi_A(\sigma) \Bigr|\\
&=&\Bigl| \sum_{A \ni 0} \Phi_A(\tau_x \sigma) \;-\sum_{A \ni 0, A
\subset (\la_n - x)} \Phi_A(\tau_x \sigma) \Bigr| \\
&\leq& \sup_{\xi} \Bigl| \sum_{A \ni 0, A \cap (\la_n - x)^c \neq
\emptyset} \Phi_A( \xi) \Bigr|.
\end{eqnarray*}
Pick $\epsilon >0$ and choose $\Delta$ such that
\begin{displaymath}
\sup_{\xi} \Bigl| \sum_{A \ni 0,A \cap \Delta^c \neq \emptyset}
\Phi_A(\xi) \Bigr|\leq \epsilon
\end{displaymath}
then
\begin{displaymath}
 \Bigl| \sum_{A \ni 0, A \cap (\la_n - x)^c \neq
\emptyset} \Phi_A( \xi) \Bigr| \; \leq \;\left\{
\begin{array}{lll} \epsilon \; \; &\rm{if}& \; \; (\la_n - x) \supset \Delta\\
\\
C \; \; &\rm{if} & \; \; (\la_n - x) \cap \Delta^c \neq \emptyset
\end{array}
\right.
\end{displaymath}
where
\begin{displaymath}
C=\sup_{\xi} \Bigl| \sum_{A \ni 0} \Phi_A(\xi) \Bigr| \; < \;
\infty.
\end{displaymath}
Since for any $\Delta \subset \mathbb{Z}^d$ finite,
\begin{displaymath}
\lim_{n \to \infty} \epsilon \; \frac{| \{x: \Delta + x \cap
\la_n^c \neq \emptyset \}|}{|\la_n|} = 0
\end{displaymath}
we obtain
\begin{displaymath}
\limsup_n \frac{1}{|\la_n|} \; \sum_{x \in \la_n} \sup_{\xi}
\Bigl| \sum_{x \ni x,A \cap \la_n^c \neq \emptyset} \Phi_A(\xi)
\Bigr| \; \leq \; \epsilon
\end{displaymath}
which by the arbitrary choice of $\epsilon>0$ proves
(\ref{chp:e4.3.1}) and the statement of the lemma. \epr

 To derive Lemma \ref{chp:lem2} from Lemma \ref{chp:lem3}, we only
have to prove that for all $\nu \in \gee_{\rm{inv}}(\gamma)$,
$e_{\nu}^+$ exists and is independent of $\gamma$. For such a
measure $\nu$, write
\begin{displaymath}
\nu(+_{\la})=\int_{\Omega} \frac{e^{-H^{\eta}_{\la_n}(+)}}{Z_{\la_n}(\eta))}\nu(d\eta)
\end{displaymath}
where $H_{\la_n}^{\eta}$ is defined \emph{via} the vacuum
potential of $\gamma$ in (\ref{chp:eqvac}). We use  Lemma
\ref{chp:lem3} to write
\begin{displaymath}
\nu(+_{\la}) \cong \int_{\Omega} \frac{e^{-H_{\la}^+(+)}}{Z_{\la}^+} \nu(d\eta)
\end{displaymath}
where $a_{\la} \cong  b_{\la}$ means $\lim_{\la} \frac{1}{|\la|} |
\log \frac{a_{\la}}{b_{\la}} |=0$. Since $\Phi$ is the vacuum
potential with vacuum state $+$, $H_{\la}^{+}(+_{\la})=0$ and
hence
\begin{displaymath}
\nu(+_{\la})=(Z_{\la}^+)^{-1}=(Z_{\la}^{\rm{free}})^{-1} =\Big[\sum_{\sigma \in
\Omega_{\la}} \exp (-\sum_{A \subset \la} \Phi_A(\sigma) ) \Big]^{-1}
\end{displaymath}
where $Z_\la^+$ (resp. $Z_\la^{\rm{free}}$) is the partition function with
$+$ (resp. free) boundary condition,
which in our case coincide. Fix $R>0$ and put
\begin{eqnarray*}
\Phi^{(R)}_A(\sigma) :&=&  \Phi_A(\sigma) \; \;  \rm{if}  \; \;
\rm{diam}(A) \leq R\\
&=& 0 \; \; \rm{if}\; \; \rm{diam}(A) \leq R
\end{eqnarray*}
then, using existence of pressure for 
finite range potentials, cf.
 \cite{I2},
\begin{displaymath}
\lim_{\la} \frac{1}{|\la|} \log Z_{\la}^{\rm{free}}(\Phi^{(R)}):=
P(\Phi^{(R)}) \; \; \rm{exists.}
\end{displaymath}
Now use

\begin{eqnarray*}
\log \frac{\sum_{\sigma} \exp{(-\sum_{A \subset \la}
\Phi_A(\sigma))}}{\sum_{\sigma} \exp{(-\sum_{A \subset \la}
\Phi_A^{(R)}(\sigma))}} \; & \leq &\; \sup_{\sigma} \Big|\sum_{A
\subset \la,\rm{diam}(A) > R} \Phi_A(\sigma) \Big|\\
& \leq & \sup_{\sigma}  \sum_{x \in \la} \Big| \sum_{A \ni x,
\rm{diam}(A)>R} \Phi_A(\sigma) \Big|\\
&\leq & \sum_{x \in \la} \sup_{\sigma} \Big|\sum_{A \ni x,
\rm{diam}(A)>R} \Phi_A(\sigma) \Big|\\
&=& |\la| \sup_{\sigma} \Big|\sum_{A \ni 0, \rm{diam}(A)>R}
\Phi_A(\sigma) \Big|
\end{eqnarray*}
and
\begin{displaymath}
\frac{\sum_{\sigma} \exp{(-\sum_{A \subset \la}
\Phi^{(R)}_A(\sigma))}}{\sum_{\sigma} \exp{(-\sum_{A \subset \la}
\Phi_A^{(R')}(\sigma))}} \leq |\la| \sup_{\sigma}\Big|\sum_{A \ni
0, \rm{diam}(A)>R \bigwedge R'} \Phi_A(\sigma) \Big|
\end{displaymath}
to conclude that $\{P(\Phi^{(R)}), R>0 \}$ is a Cauchy net
with limit
\begin{displaymath}
\lim_{R \to \infty} P(\Phi^{(R)})=\lim_{\la \uparrow \mathbb{Z}^d}
\frac{1}{|\la|} \log Z_{\la}^{\rm{free}}=e_{\nu}^+
\end{displaymath}
which depends only on the vacuum potential (hence on the
specification $\gamma$). This proves that $e_{\nu}^+$ and
$e_{\mu}^+$ exist for all $\mu, \nu \in \gee_{\rm{inv}}(\gamma)$,
and depends of $\gamma$ only. Therefore,
\begin{displaymath}
\lim_{\la \uparrow \Z^d} \frac{1}{|\la|} \log
\frac{\mu(+_{\la})}{\nu(+_{\la})}=e_{\nu}^+-e_{\mu}^+=0
\end{displaymath}
which proves Lemma \ref{chp:lem2}.
\epr

A direct consequence of this lemma is that in the framework
of Theorem \ref{chp:thm7}, $e_{\nu}^+$ exists and conditions C1
and C2 are true. We obtain the theorem by applying Theorem
\ref{chp:thm4}.
 \epr

\section{Examples}
\subsection{The GriSing random field}
The GriSing random field is an example of joint measure of disordered systems, 
studied in more extent  in Section 6. It has been studied in \cite{EMSS} and provides 
 an easy example of a non-Gibbsian
random fields which fits in the framework of our theorems. The random
field is constructed as follows.
Sites
are empty or occupied according a Bernoulli product measure of
parameter $p<p_c$ where $p_c$ is the percolation
threshold for site percolation on $\Z^d$. For any realization
$\eta$ of occupancies where all occupied clusters are finite,
we have the Gibbs measure on configurations $\si\in \{-1 ,+1\}^{\Z^d}$
\[
\mu^\eta_\beta (d\si)
\]
which is the product of free boundary condition Ising measures on
the occupied clusters. More precisely, under $\mu^\eta_\beta$ spin configurations
of occupied clusters $C$ are independent and distributed as:
\[
\mu_{\beta, C} (\si_C) = \frac{1}{Z_{\la}} \; e^{-\beta \sum_{\langle xy
\rangle \subset C} \sigma(x) \sigma(y) }
\]
The GriSing random field is then defined as:
\[
\xi (x) = \si (x) \eta (x).
\]
In words, $\xi(x)=0$ for unoccupied sites and equal to the spin
$\si (x)$ at occupied sites.

We denote by $K_{p,\beta}$ the law of the random field
$\xi$.

It is known  that for any $p\in (0,1)$,
$\beta$ large
enough, $K_{p,\beta}$ is not a Gibbs measure
(see \cite{EMSS} for $p<p_c$ and \cite{ku} for any $p\in (0,1)$). The points
of essential discontinuity of the conditional probabilities
$K_{p,\beta}(\si(0)|\xi_{\Z^d\setminus\{0\}})$ are a subset
of
\[
D= \{ \xi: \xi\ \mbox{contains an infinite cluster of
occupied sites}\}.
\]
Since $p<p_c$,
there exists a specification $\gamma$ such that $\{K_{p,\beta}\}=\geeg$
and such that for the continuity points $\Omega_\gamma$, we have
$K_{p,\beta}(\Omega_\gamma) =1$, i.e., $K_{p,\beta}$ is almost
Gibbs. Moreover, if we choose $\xi_0 \equiv 0$ as a telescoping reference
configuration, then clearly $\si\in D^c$ implies
$\si^{\xi_0}\in D^c$, i.e., in this case
$\Omega_\gamma\subset \Omega^{<0}_\gamma$.
 Therefore, in this example condition C1 is satisfied as soon
as $\mu$ concentrates on $D^c$. Using $\{K_{p,\beta}\}=\geeg$, and

\[
\lim_{\la\uparrow\Z^d}\frac{1}{|\la|}\log K_{p,\beta}(0_\la)=\log (1-p)
\]
we obtain the following proposition:

\bp
If $\mu(D) =0$ then $h(\mu|K_{p,\beta})$ exists
and is zero if and only if $\mu=K_{p, \beta}$.
\ep

\subsection{Decimation}
Let $\mu^+_\beta$ (resp. $\mu^-_\beta$)
be the  low-temperature ($\beta>\beta_c$) plus (resp. minus)
phase of the Ising model on $\Z^d$. For $b\in\N$, $\nu^+_\beta$ denotes
its decimation, i.e., the distribution of $\{ \si( bx): x\in\Z^d\}$ when
$\si$ is distributed according to $\mu^+_\beta$. It is known
that $\nu^+_\beta$ is not a Gibbs measure \cite{EFS}. In
\cite{FP} it is proved that there exists a monotone
specification $\gamma^+$ such that $\nu^+_\beta\in\gee(\gamma^+)$.
In \cite{FLNR} it is proved that the points of continuity
$\Omega_{\gamma^+}$ satisfy $\nu^+_\beta(\Omega_{\gamma^+})=1$, i.e., $\nu^+_\beta$ is almost Gibbs.
The point of continuity of $\gamma^+$ can be described as those
configurations $\eta$ for which the ``internal spins" do not exhibit
a phase transition when the decimated spins are fixed to
be $\eta$. E.g., the all plus and the all minus configurations
are elements of $\Omega_{\gamma^+}$, but the alternating configuration
is not.

The first part of the
variational principle for $(\gamma^+,\nu^+_\beta,\mathcal{M})$ has already been proved in \cite{FLNR} (and is direct by Theorem \ref{chp:thm1}), with  a set $\mathcal{M}$  consisting of  the translation invariant measures which concentrate on $\Om_{\gamma^+}$ 
. Here we complete this result by adding a second part:
\begin{theorem}\label{chp:thm10}
For any $\mu \in
\minv$ satisfying C1 for $\gamma^+$,
\begin{enumerate}
\item $h(\mu|\nu)$ exists.
\item We have the equivalence
\begin{displaymath}
\mu \in \gee_{\rm{inv}}(\gamma^+) \; \; \Longleftrightarrow \; \;
h(\mu|\nu^+)=0.
\end{displaymath}
\end{enumerate}
\end{theorem}

We first use a lemma.
\begin{lemma}
$\mu \in \gee(\gamma^+)$ and $\mu(\Om_{\gamma^+})=1$ implies
\be \label{chp:fkg}
\nu_{\beta}^+ \; \preceq \; \mu \; \preceq \; \nu_{\beta}^+.
\ee
\end{lemma}
\bpr
Consider $f$ monotone. By monotonicity of $\gamma^+$ \cite{FP}, for all $\la \in \s$,
\[ \int f d\mu = \int_{\Om} (\gamma_{\la}^+ f)(\omega) \mu(d \om) \leq  \int_{\Om} (\gamma_{\la}^+ f)(+) \mu(d \om)=(\gamma_{\la}^+ f)(+).
\]
Taking the limit $\la \uparrow \Z^d$ gives, and using $\gamma_{\la}^+(\cdot|+)$ goes to $\nu_{\beta}^+$,
\[
\int f d\mu \leq \int f d\nu_{\beta}^+.
\]
Similarly, using $\mu(\Om_{\gamma})=1$, and the expression of $\Om_\g$ in \cite{FP}, we have $\gamma^+(f)=\gamma^-(f)$, $\mu$-a.s. and hence
\[
\int f d\mu =\int \gamma^-_\la (f) d\mu \; \geq \; \gamma^-_\la f(-)
\]
which gives
\[
\int f d\mu \geq \int f d\nu_{\beta}^-.
\]
\epr

The following corollary proves Theorem \ref{chp:thm10} using Theorem \ref{chp:thm4}.
\bp
\ben
\item $e_{\nu_{\beta}}^+=-\lim_{\la \uparrow \Z^d} \frac{1}{|\la|} \log \nu_{\beta}^+(+_\la)$ exists.
\item For any $\mu \in \gee(\gamma^+)$,
\[
\lim_{\la \uparrow \Z^d} \frac{1}{|\la|} \log \frac{\mu(+_\la)}{\nu_{\beta}^+(+_\la)}=0.
\]
\een
\ep
\bpr
\ben
\item Follows from subadditivity and positive correlations.
\item Follows from stochastic domination (\ref{chp:fkg}) and
\[
\lim_{\la \uparrow \Z^d} \frac{1}{|\la|} \log \frac{\nu^+(+_\la)}{\nu_{\beta}^-(+_\la)}=\lim_{\la \uparrow \Z^d} \frac{1}{|\la|} \log \frac{\mu_{\beta}^+(+_{b \la})}{\mu_{\beta}^-(+_{b \la})}=0
\]
where, to obtain the last equality,
we used that $\mu^+_\beta, \mu^-_\beta$ are the Ising
plus and minus phases.
\een
\epr

\br

We conjecture that $C1$ is satisfied for
any ergodic measure
$\mu\in\gee (\gamma^+ )$ in dimension $d=2$. This amounts of to prove that
the internal spins do not show a phase transition, given
a ``typical configuration of $\mu$" on $b\Z^d$
to the left of the origin,
and all $+$ on $b\Z^d$
to the right. Fixing these decimated spins acts
as a magnetic field, pushing the spins on
the right of the origin into a ``plus-like" phase and
the spins on the left of the origin in a ``plus-like"
or ``minus-like" phase, depending on $\mu$. The location
of the interface between ``right and left" should not
depend on the boundary condition in $d=2$ (no Basuev transition).
However, we do not have
a rigorous proof of this fact.
\er

\section{Examples II: Joint measure of random spin systems}

We consider  the joint measures of disordered spins-systems on the
product of spin-space and disorder-space defined in terms of a
quenched absolutely convergent Gibbs-interaction and an a
priori-distribution of the disorder variables. They were treated
before in \cite{ku,ku1} and provide a broad class of examples of
generalized Gibbs measures. A specific example of this, the
GriSing field,  was already considered in Section 5.1.

First we prove that, for the same quenched potential,
the relative entropy density between corresponding,
possibly different
joint measures is always zero.
Next we prove in generality
that these measures are asymptotically decoupling
whenever the a-priori distribution of the disorder is.
The useful notion of
asymptotic decoupled was recently coined by Pfister (2002), and provides a broad class of measures, including local transformations of Gibbs measures, for which the existence of relative entropy density and the large deviation principle holds.
Using results these results, we easily obtain
existence of the relative entropy density.
Next we specialize to the specific example
of the random field Ising
model in Section 6.3. We focus on the interesting region of the parameter
space when there is a phase transition for
the spin-variables, for almost any configuration
of disorder variables.
Here we show on the basis of \cite{ku1}  that the
joint plus and the joint minus state for the
same quenched potential are not compatible with the same interaction
potential.
In \cite{ku1} it was already shown that there is always
a translation-invariant convergent potential, or a
possibly non-translation-invariant absolutely convergent
potential for the corresponding joint measure.
We also discuss this in more detail and
sketch a proof on the basis of \cite{ku1} and
the RG-analysis
of Bricmont and Kupiainen (1988)  that shows that there is a translation-invariant
joint potential that
even decays like a stretched exponential.
This provides an explicit example of a weakly (but not almost)
Gibbsian measure for which the variational principle fails.

\subsection{Setup}

\vskip 0.1in

We consider disordered models of the following
general type.
We assume that the configuration space of the quenched
model is again as detailed in Section 2.1 and we denote
the spin variables by $\si$. Additionally
we assume that there are also disorder
variables $\eta=(\eta_x)_{x\in \Z^d}$ entering the game,
taking values in an infinite
product space $(E')^{\Z^d}$, where again $E'$ is
a finite set.
We denote the {\it joint variables} by $\x=(\x_x)_{x\in \Z^d}
=(\si,\eta)=(\si_x ,\eta_x)_{x\in \Z^d}$.
It will be convenient later also to write simply $(\si\eta)$
to denote the pair $(\si,\eta)$.

One essential ingredient of the model is given
by the {\it defining potential}
$\Phi=(\Phi_A)_{A\sb \Z^d}$ depending on the joint variables
$\x=(\si,\eta)$. $\Phi_A(\x)$ depends on $\x$ only through $\x_{A}$.
We assume that $\Phi$ is finite range.
When we fix a realization of the disorder $\eta$,
we have a potential for the spin-variables $\si$ that is typically
non-translation invariant. We then define the corresponding
{\it quenched Gibbs specification} by Definition \ref{chp:def3}
using the notation
\be\label{chp:6.1.2}
\mu_{\L}^{\bar \si}[\eta](B)
:=\frac{1}{Z^{\bar \si}_{\L}[\eta]}
\sum_{\si_\L}1_{B}(\si_{\L}{\bar \si}_{\Z^d\ba \L})
e^{-\sum_{A:A\cap \L\neq \emp}
\Phi_{A}( \si_{\L} \bar \si_{\Z^d\ba \L},
\eta)}.
\ee
Here we do not make the defining potential explicit
anymore in order not to overburden notion.
The measures (\ref{chp:6.1.2}) are also
called more loosely {\em quenched finite-volume Gibbs measures}.
Obviously, the finite-volume
summation is over $\si_{\L}\in E^{\L}$.

The second ingredient of the quenched model
is the distribution of the disorder
variables $\P(d\eta)$.
Most of the times  in the theory of disordered systems
one considers the case of i.i.d. variables, but we can and
will be more general here.

The objects of interest will then be the infinite volume
{\it joint measures} $K^{\bar \si}(d\x)$,
by which we understand
any limiting measure
of $\lim_{\L\uparrow\Z^d}\P(d\eta)\mu_{\L}^{\bar \si}[\eta](d\si)$ in the product
topology on the space of joint variables.
Of course, there are examples for different joint measures
of the same quenched Gibbs specification
for different spin boundary conditions $\bar \si$.
In principle there can even be different ones for
the same spin-boundary condition $\bar \si$, depending
on the sub-sequence.

For all of this the reader might think  of the concrete
example of the {\it random field Ising model.}
Here the spin variables $\si_x$ take values in
$\{-1,1\}$.
The disorder variables are given by the random fields $\eta_x$
that are i.i.d. with single-site distribution $\P_0$
that is supported on a finite set $\HH_0$ and assumed
to be symmetric.
The defining potential $\Phi(\si,\eta)$ is
given by $\Phi_{\{x,y\}}(\si,\eta)=-\b\si_x\si_y$
for nearest neighbors $x,y\in \Z^d$,
$\Phi_{\{x\}}(\si,\eta)=-h\eta_x\si_x$,  and $\Phi_A=0$ else.

\subsection{Relative entropy for joint measures}
\medskip

For the first result we do not need
the independence of the disorder field. In fact,
without any decoupling assumption on $\P$ we have the following.

\begin{theorem}\label{chp:thm6.2}
Denote by $K^{\bar\si}$ and $K^{\bar\si'}$ two
joint measures
for the same quenched Gibbs specification $\mu^{\Large{\bf{\cdot}}}_{\L}[\eta](d\si)$,
obtained with any two spin boundary conditions $\bar \si$
(and $\bar \si'$ respectively),
along any subsequences $\L_N$ (and $\L'_N$ respectively).
Then their relative entropy density vanishes, i.e.,
$h(K^{\bar\si}|K^{\bar\si'})=0$.
\end{theorem}

\br
Note that we are more general than in the usual set up and we do not need to assume translation invariance, not even of the defining potential $\Phi$.
\er

\br
This result is neither directly related to the first part nor to the second part of the variational principle. It does not yield the first part (which will be proved differently)
because it is not clear that every measure that is compatible with the same specification as $K^{\bar\si'}$ can be written in terms of $K^{\bar\si}$. Applied to the random field Ising model in Section 6.3, this result will disprove the second part of the variational principle for weakly but not almost Gibbs measures.
\er

\bpr
We have from the definition of the joint measures as limit points
with suitable sequences of volumes
\be \label{chp:6.1.3}
\frac{K^{\bar\si}(\si_{\L}\eta_{\L})}{K^{\bar\si'}(\si_{\L}\eta_{\L})}
=\frac{\lim_{N}K^{\bar\si}_{\L_N}(\si_{\L}\eta_{\L})}{\lim_{N}K^{\bar\si'}_{\L'_N}(\si_{\L}\eta_{\L})}
=\frac{\lim_{N}\int\P(d\tilde\eta)1_{\eta_{\L}}\mu_{\L_N}^{\bar \si}[\tilde\eta](\si_{\L})
}
{\lim_{N}\int\P(d\tilde\eta)1_{\eta_{\L}}\mu_{\L_N}^{\bar \si'}[\tilde\eta](\si_{\L})}
\ee
Here and later we will write in short
$1_{\eta_{\L}}$ for the indicator function
of the event that the integration variable
$\tilde \eta$ coincides with the fixed configuration $\eta$
on $\L$. We have from the finite range of the disordered potential that
\[ 
\sup_{\si\eta=\si'\eta'\hbox { on }\L}\Bigl|\sum_{A}\bigl(
\Phi_A(\si\eta)-\Phi_{A}(\si'\eta')\bigr)\Bigr|\leq C_1 |\partial \L|
\]
for cubes $\L$ with some finite constant $C_1$. By $\del \la$ we mean the $r$-boundary of $\la$, where $r$ is the range of $\Phi$.
So we get that for $N$ large enough
\[
e^{ -2 C_1 |\partial \L|}
\mu_{\L}^{\hat \si}[\eta_{\L}\hat\eta_{\Z^d\ba\L}](\si_{\L})\leq
\mu_{\L_N}^{\bar \si}[\eta_{\L}\tilde\eta_{\Z^d\ba\L}](\si_{\L})\leq
e^{ 2 C_1 |\partial \L|}
\mu_{\L}^{\hat \si}[\eta_{\L}\hat\eta_{\Z^d\ba\L}](\si_{\L}) 
\]
for any joint reference configuration $\hat\si\hat\eta$.
This gives the upper bound $ e^{ 4 C_1 |\partial \L|}$
on the r.h.s. of (\ref{chp:6.1.3}), by application of the last
inequalities on numerator and denominator of (\ref{chp:6.1.3})
for the same reference configuration.

This implies for the finite-volume relative entropy
an upper bound of the order of the boundary, i.e.,

\[
h_{\L}(K^{\bar\si}|K^{\bar\si'})=\sum_{\si_{\L}\eta_{\L}}K^{\bar\si}
(\si_{\L}\eta_{\L})
\log\frac{K^{\bar\si}(\si_{\L}\eta_{\L})}{K^{\bar\si'}(\si_{\L}\eta_{\L})}
\leq  4 C_1 |\partial \L|.
\]

From that clearly follows the claim
$h(K^{\bar\si}|K^{\bar\si'})\leq \limsup_{n\uparrow \infty}\frac{1}{|\L_n|}
h_{\L_n}(K^{\bar\si}|K^{\bar\si'})=0$ for $(\L_n)_{n \in \N}$ a sequence of cubes.
\epr

Also the next theorem can be proved in a natural way when
we relax the independence assumption of
the a priori distribution $\P$ of the
disorder variables. It says that
the property of being {\it asymptotically decoupled}
carries over from the distribution of the
disorder fields to any corresponding
joint distribution. Following \cite{P}, we give the following

\bd\label{AD} A probability measure $\P\in\minv$ is called
asymptotically decoupled (AD) if there exists sequences $g_n$,
$c_n$ such that
\[
\lim_{n\to\infty} \frac{c_n}{|\la_n|} =0,\; \;  \lim_{n\to\infty}
\frac{g_n}{n} =0
\]
and  for all $A\in \fe_{\la_n}$, $B\in \fe_{\la_{n+g_n}^c}$ with
$\P (A)\P (B) \not= 0$: \be\label{addef} e^{-c_n} \leq
\frac{\P (A \cap B)}{\P (A) \P (B)} \leq e^{c_n}. \ee \ed

\begin{theorem}\label{chp:Theorem 6.8}Suppose $\P$ is asymptotically decoupled with functions $g_n$ and $c_n$.
Assume that $K^{\bar \si}$ is a corresponding translation
invariant joint measure of a quenched random system, with a defining finite range potential.
Then $K^{\bar \si}$ is asymptotically decoupled
with functions $g'_n=g_n$ and $c'_n=c_n+ C |\partial \L_n|$, where $C$ is a real constant.
\end{theorem}
\bpr It suffices to show
that for any {\em finite} $V\subset \L_{n+g'(n)}^c$ we have
\be \label{2.1}
e^{-c'_n}\leq\frac{K(\x_{\L_n}\x_{V})}{K(\x_{\L_n}) K(\x_{V})}=\frac{K(\si_{\L_n}\eta_{\L_n}\si_{V}\eta_{V})}{K(\si_{\L_n}\eta_{\L_n}) K(\si_{V}\eta_{V})}
\leq e^{c'_n}.
\ee
We only show the upper bound.
It suffices to show
\[
\limsup_{n}
\frac{K^{\bar \si}_{\tilde\L_n}(\si_{\L_n}\eta_{\L_n}\si_{V}\eta_{V})}{K^{\bar \si}_{\tilde\L_N}(\si_{\L_n}\eta_{\L_n})K^{\bar \si}_{\tilde\L_N}(\si_{V}\eta_{V})}\leq e^{c_n}
\]
for any sequence $\tilde \L_N$. The quantity under the limsup equals
\be \label{ad}
\frac{\int \P(d\tilde \eta) 1_{\eta_{\L_n}}1_{\eta_{V}}
\mu_{\tilde \L_N}^{\bar \si}[\tilde\eta](\si_{\L_n}\si_V)}
{
\int \P(d\tilde \eta_1) 1_{\eta_{\L_n}}
\mu_{\tilde \L_N}^{\bar \si}[\tilde\eta_1](\si_{\L_n})
\int \P(d\tilde \eta_2) 1_{\eta_{V}}
\mu_{\tilde \L_N}^{\bar \si}[\tilde\eta_2](\si_V)}.
\ee
Look at the term under the disorder-integral in the numerator.
We have by the compatibility of the quenched kernels that

\begin{eqnarray*}
& &\mu_{\tilde \L_N}^{\bar \si}[\eta_{\L_n}\eta_{V}\tilde\eta_{\Z^d\ba(\L_n\cup V)}]
(1_{\si_{\L_n}}1_{\si_V})\\
&=& \int \mu_{\tilde \L_N}^{\bar \si}[\eta_{\L_n}\eta_{V}\tilde\eta_{\Z^d\ba(\L_n\cup V)}](d\tilde \si)
1_{\si_V}
\mu_{\L_n}^{\tilde\si}[\eta_{\L_n}\eta_{V}\tilde\eta_{\Z^d\ba(\L_n\cup V)}]
(1_{\si_{\L_n}})\\
&\leq &
e^{ 2 C_1 |\partial \L_n|}
\mu_{\L_n}^{\hat \si}[\eta_{\L_n}\hat\eta_{\Z^d\ba\L_n}](\si_{\L_n})
\times \mu_{\tilde \L_N}^{\bar \si}[\eta_{\L_n}\eta_{V}\tilde\eta_{\Z^d\ba(\L_n\cup V)}](1_{\si_V})\\
\end{eqnarray*}
where the inequality follows from the uniform absolute
convergence of the quenched potential, for any reference
configuration $\hat \si \hat\eta$.

We use that
\[
\mu_{\tilde \L_n}^{\bar \si}[\eta_{\L_n}(\tilde\eta_1)_{\Z^d\ba \L_n}](\si_{\L_n})
\geq e^{ -2 C_1 |\partial \L_n|}\mu_{\L_n}^{\hat \si}[\eta_{\L_n}\hat\eta_{\Z^d\ba\L_n}](\si_{\L_n})
\]
and the similar lower bound
on the first disorder-integral
in the denominator of (\ref{ad})  with the same reference
joint reference configuration $\hat\si\hat\eta$. From this
we get an upper bound on (\ref{ad}) in the form of

\be \label{2.6}
e^{ 4 C_1 |\partial \L_n|}
\frac{\int \P(d\tilde \eta) 1_{\eta_{\L_n}}1_{\eta_{V}}
\mu_{\tilde \L_N}^{\bar \si}[\tilde\eta](\si_V)}
{\int \P(d\tilde \eta_1) 1_{\eta_{\L_n}}
\int \P(d\tilde \eta_2) 1_{\eta_{V}}
\mu_{\tilde \L_N}^{\bar \si}[\tilde\eta_2](\si_V)}.
\ee

Last we need to control the influence of the
variation of the random fields inside the finite volume
$\eta_{\L_n}$ on the Gibbs-expectation outside.
We have that
\[
\mu_{\tilde \L_N}^{\bar \si}[\eta_{\L_n}\tilde\eta_{\Z^d\ba \L_n}](\si_V)
\leq e^{2 C_1 |\partial \L_n|}
\mu_{\tilde \L_N}^{\bar \si}[\eta_{\L_n}^{(1)}\tilde\eta_{\Z^d\ba \L_n}](\si_V)
\]
for any configurations $\eta$ and $\eta^{(1)}$ inside $\L_n$.
This gives the following upper bound on (\ref{2.6})
\[
e^{ 8 C_1 |\partial \L_n|}
\frac{\int \P(d\tilde \eta) 1_{\eta_{\L_n}}1_{\eta_{V}}}
{\int \P(d\tilde \eta_1) 1_{\eta_{\L_n}}
\int \P(d\tilde \eta_2) 1_{\eta_{V}}
}.
\]
But this, by the property of asymptotic decoupling of
the disorder field, is bounded by
$e^{ 8 C_1 |\partial \L_n|+c_n}$
and the proof of the upper bound in (\ref{2.1}) is done.
The proof of the lower bound  is similar.
\epr

\bigskip

Applying Pfister's theory \cite{P}, we have
\bigskip
\bc Suppose $\P$ is asymptotically decoupled and that $K^{\bar \si}$ is a corresponding translation
invariant joint measure of a quenched random system, with a defining finite range potential. Then
$h(K|K^{\bar \si})$ exists
for all translation invariant probability measures $K$.
\ec
Moreover we have the following explicit formula:
\begin{theorem}\label{chp:6.3}
Suppose that the defining potential $\Phi(\si,\eta)$ is translation
invariant and that $\P$ is asymptotically decoupled.
Suppose that $K^{\bar \si}$ is translation invariant joint
measure constructed with the boundary condition $\bar \si$.
Suppose that $K$ is a translation-invariant
measure on the product space. Denote by $K_d$ its marginal
on the disorder variables $\eta$.
Then
\begin{eqnarray*}
 h(K|K^{\bar\si})&=& h(K_d|\P) - h(K)-h(K_d)\\
&+&\sum_{A\ni 0}\frac{1}{|A|}K\Bigl(\Phi_A(\si\,\eta=\cdot)
\Bigr)
+ K\Bigl(
\lim_{\L}\frac{1}{|\L|}
\log Z^{\bar \si}_{\L}(\eta=\cdot) \nonumber
\Bigr)
\end{eqnarray*}
where $h(K)$
is the Kolmogorov-Sinai entropy (\ref{chp:KS}).
\end{theorem}

\br
The third term has the meaning of the $K$-expectation
of the `joint energy'.
The last term is the $K$-mean of the ``quenched pressure''.
Note that it is boundary condition $\bar \si$-independent, of course.
 \er

\br
In the case that $\P$ is
a Gibbs distribution, the existence
of the relative entropy density is obtained
directly, i.e., without relying on Pfister's theory.
 \er

\bpr We have
\[
\frac{1}{|\L|}h_{\L}(K|K^{\bar\si})
=\frac{1}{|\L|}
\sum_{\si_{\L} \eta_{\L}} K(\si_{\L} \eta_{\L})
\log K(\si_{\L} \eta_{\L})
-\frac{1}{|\L|}
\sum_{\si_{\L} \eta_{\L}} K(\si_{\L} \eta_{\L})
\log K^{\bar \si}(\si_{\L} \eta_{\L})
\]
where the first term converges to $-h(K)$.
For the second term we use the approximation
\[
\sup_{\bar \si,\hat\si,\hat\eta}
\Biggl|\log\Biggl( \frac{K^{\bar \si}(\si_{\L} \eta_{\L})}
{\P(\eta_{\L})\mu_{\L}^{\hat \si}[\eta_{\L}\hat\eta_{\Z^d\ba\L}](\si_{\L})}
\Biggr)\Biggr|
\leq  2 C_1 |\partial \L|.
\]
First we have
\[
-\frac{1}{|\L|}
\sum_{\si_{\L} \eta_{\L}} K(\si_{\L} \eta_{\L})
\log \P(\eta_{\L})=\frac{1}{|\L|}h_{\L}(K_d|\P)
- \frac{1}{|\L|}
\sum_{\eta_{\L}} K_d(\eta_{\L})\log K_d(\eta_{\L}).
\]
The second term converges to $h(K_d)$. The first
term converges to $h(K_d|\P)$. This is clear either by
the classical
theory for the case that $\P$ is Gibbs or even independent,
or by Pfister's theory if $\P$ is asymptotically decoupled.
Next, by definition
\[
\log\mu_{\L}^{\hat \si}[\eta_{\L}\hat\eta_{\Z^d\ba\L}](\si_{\L})
=-\sum_{A:A\cap \L\neq \emp}
\Phi_{A}( \si_{\L}\hat \si_{\Z^d\ba \L}\eta_{\L}\hat\eta_{\Z^d\ba\L})
-\log Z^{\hat \si}_{\L}(\eta_{\L}\hat\eta_{\Z^d\ba \L}).
\]
Using translation-invariance of the measure $K$ we get
that the application of  $\frac{1}{|\L|}\int K(d\si_{\L}
d\eta_{\L})$ over the first sum of the r.h.s.
converges to
$-\sum_{A\ni 0}\frac{1}{|A|} K(
\Phi_A(\si \, \eta=\cdot))$.
To see that the average over the last term converges
we use the ergodic decomposition
of $K_d$ to write $K_d(d\eta)=\int\r(d\k) \k(d\eta)$ where
$\r(d\k)$ is a probability measure
that is concentrated on the ergodic measures on $\eta$.
Fix any ergodic $\k$. For
$\k$-a.e. disorder configuration $\eta$
we have the existence of the limit
$-\lim_{\L}\frac{1}{|\L|}
\log Z^{\bar \si}_{\L}(\eta=\cdot)$,
by standard arguments \cite{sepp}.
The convergence is also in $L^1$, by dominated
convergence. So we may integrate over $\r$
to see the statement of the theorem.
\epr

\subsection{Discussion of the first part of the variational principle for joint measures}

To discuss the first part of the variational principle we will use an
explicit representation of the conditional expectations of the joined measures.
For this we need to restrict to the case that $\P$
is a product measure.
First, in the situation detailed below, we prove the first part
of the variational principle by direct arguments.
Next, we illustrate the criteria given in the general
theory of Section 3.4 by showing that they can be
verified in the context of joint measures in the almost Gibbsian case,
giving then an alternative proof of the variational principle.

We start with the following proposition of \cite{ku1}. \\

\bp\label{kulsprop} 
 Assume that $\P$
is a product measure. Assume that
there is a set of realizations of $\eta$'s
of $\P$-measure one
such that the quenched infinite-volume Gibbs measure
$\mu[\eta]$ is a weak limit
of the quenched finite-volume measures (\ref{chp:6.1.2}).
Then, a version of the infinite-volume
conditional expectation of the
corresponding joint measure $K^{\mu}(d\si,d\eta)=\P(d\eta)\mu[\eta](d\si)$
is given by the formula
\be \label{chp:6.4.2}
K^{\mu}\left[\x_{\la}\bigl|\x_{\la^c}\right]
= \frac{\mu^{\hbox{\rm ann,}\x_{\del \la}}_{\la}(\x_{\la})}{
\int\mu^{\hbox{\rm ann,}\x_{\del \la}}_{\la}(d\tilde\eta_{\la})
Q^\mu_{\la}(\eta_{\la},\tilde\eta_{\la},\eta_{\la^c})}.
\ee
Here $\mu^{\hbox{\rm ann,}\x_{\del \la}}_{\la}(\x_{\la})$
is the {\em trivial annealed local specification} given
by in terms of the potential $U^{\rm{triv}}_{A}(\si,\eta)=
\Phi_{A}(\si,\eta)-1_{A=\{x\}}\log \P_0(\eta_x)$ w.r.t counting
measure on the product space.

Further we have put
\[
Q^{\mu}_{\la}(\eta^{1}_{\la},\eta^{2}_{\la},\eta_{\la^c})
=\mu[\eta^2_{\la}\eta_{\la^c}]
(e^{-\D H_{\la}(\eta^1_{\la},\eta^2_{\la},\eta_{\del {\la}})})
\]
 where
 \[\Delta H_{\la}(\eta^{1}_{\la},\eta^{2}_{\la},\eta_{\la^c})(\si)=\sum_{A \cap \la \neq \emp}\Bigl (\Phi_A(\si,\eta_{\la}^1 \eta_{\la^c} ) - \Phi_A(\si,\eta_{\la}^2 \eta_{\la^c}) \Bigr).
\]
\ep
According to our assumption on the measurability
on $\mu[\eta]$, $Q^{\mu}_{\la}$ depends measurably on $\eta_{\la^c}$.
We fix a version of the map  and
define the r.h.s. of (\ref{chp:6.4.2}) to
be the specification $\g^{\mu}$.
Note that for the random field Ising model, this specification
 exists for all random field configurations
by monotonicity.

In this context we always have the first part of the variational principle. Note
that we do not need any further assumption about almost Gibbsianness.

\begin{theorem} Assume that $\pee$ is a product measure.
There exists a constant $C$ depending only on $\Phi$, $\pee$ such that for
any $K,K'\in\geegu
$ one has
\[
\sup_\xi \left | \log \frac{K(\xi_\la)}{K'(\xi_\la)} \right |\leq C |\partial \la|.
\]
In particular $h(K|K')=h(K'|K)=0$.
\end{theorem}
\bpr
Using $K,K'\in\geegu$, it suffices to show that we have the estimate
\[
\frac{\gamma^\mu_\la (\xi_\la|\xi_{\la^c})}{\gamma^\mu_\la (\xi_\la|\xi'_{\la^c})}
\leq  e^{C|\partial\la|}
\]
where the constant $C$ is independent of $\la,\xi,\xi'$.
From the explicit representation (\ref{chp:6.4.2}) we obtain
\beq\label{propo1}
 \frac{\gamma^\mu_\la (\xi_\la|\xi_{\la^c})}{\gamma^\mu_\la (\xi_\la|\xi'_{\la^c})}
&=&
\frac{\mu_\Lambda^{\mbox{ann},\xi_{\partial\Lambda}} (\xi_{\Lambda})}
{\mu_\Lambda^{\mbox{ann},\xi'_{\partial\Lambda}} (\xi_{\Lambda})}
\frac{\int\mu_\Lambda^{\mbox{ann},\xi'_{\partial\Lambda}} (d\tilde{\eta}_{\Lambda})Q^\mu_\la (\eta_\la,\tilde{\eta}_\la,
\eta'_{\la^c})}
{\int\mu_\Lambda^{\mbox{ann},\xi_{\partial\Lambda}} (d\tilde{\eta}_{\Lambda})Q^\mu_\la (\eta_\la,\tilde{\eta}_\la,
\eta_{\la^c})}.
\eeq
Using the definition of $\mu_\Lambda^{\mbox{ann},\xi_{\partial\Lambda}}$ and using
the finite range assumption on $\Phi$, we obtain the bound
$e^{c |\partial\la|}$ for the first factor on the r.h.s. of (\ref{propo1}).
The second factor on the r.h.s. of (\ref{propo1}) is bounded by
\[
\left(\sup_{\tilde{\eta}_\la}\frac{Q^\mu_\la (\eta_\la,\tilde{\eta}_\la,
\eta'_{\la^c})}
{Q^\mu_\la (\eta_\la,\tilde{\eta}_\la,
\eta_{\la^c})}\right)
\frac{\int\mu_\Lambda^{\mbox{ann},\xi'_{\partial\Lambda}} (d\tilde{\eta}_{\Lambda})Q^\mu_\la (\eta_\la,\tilde{\eta}_\la,
\eta_{\la^c})}
{\int\mu_\Lambda^{\mbox{ann},\xi_{\partial\Lambda}} (d\tilde{\eta}_{\Lambda})Q^\mu_\la (\eta_\la,\tilde{\eta}_\la,
\eta_{\la^c})}.
\]
Using the same argument on
$\mu_\Lambda^{\mbox{ann},\xi_{\partial\Lambda}}$ again, we see that the second
factor is bounded by $e^{C|\partial\la|}$.
To estimate the first factor, remind the explicit expression
\begin{eqnarray*}
Q^\mu_\la (\eta_\la,\tilde{\eta}_\la,
\eta_{\la^c})
&=& \mu[\tilde{\eta}_\la\eta_{\la^c}] \left( e^{-\Delta H_\la (\eta_\la,\tilde{\eta}_\la,\eta_{\la^c})}\right)
\\
&\leq & e^{c|\partial\la|}\mu[\tilde{\eta}_\la\eta_{\la^c}]\left(e^{-\Delta H_\la (\eta_\la,\tilde{\eta}_\la,\eta'_{\la^c})}\right).
\end{eqnarray*}
Here the inequality follows from the definition of $H_\la$ and the finite range
property of $\Phi$.
Now use the definition of the quenched kernels and once again the finite range of $\Phi$
to see that the last expectation is bounded from above by
\[
e^{c|\partial\la|}\mu[\tilde{\eta}_\la\eta'_{\la^c}]\left(e^{-\Delta H_\la (\eta_\la,\tilde{\eta}_\la,\eta'_{\la^c})}\right)
= Q^\mu_\la (\eta_\la,\tilde{\eta}_\la,
\eta'_{\la^c}).
\]
This finishes the proof.
\epr
\bigskip

Let us now check what can be said about the criteria
for the first part of the variational principle for joint measures.
It turns out that it is natural to use the criteria
given in Section 3.4
with a measure $\l$ that is not a Dirac measure.
Instead, let us
take any translation invariant configuration
$\si^0$ and put $\l:=\P\otimes \d_{\si^0}$.

First, using the arguments given in the proof of Theorem \ref{chp:6.3},
it is simple in this situation to see that  the limit (\ref{elmu}) exists and to give an explicit expression for it.
\bp
Suppose that the defining potential $\Phi$ is translation
invariant. Suppose that $K^{\bar \si}$ is translation invariant joint
measure constructed with the boundary condition $\bar \si$. Then
\[
e^{\l}_{K^{\bar \si}}
= - h(\P)
+ \sum_{A\ni 0}\int\P(d\eta)\frac{\Phi_A(\si^0,\eta)}{|A|}
+ \int\P(d\eta)\lim_{\L\uparrow\Z^d}\frac{1}{|\L|}\log Z^{\bar \si}_{\L}[\eta]
\]
exists.
\ep

Put
\[
\HH_{\mu}:=\{\eta \in \HH, \eta\mapsto
Q^{\mu}_{x}(\eta^{1}_{x},\eta^{2}_{x},\eta_{\Z^d\ba x})\hbox{ is continuous }\forall
x,\eta^{1}_{x},\eta^{2}_{x}\},
\]
then we have that $\si\eta\in \O_{\g^\mu}\Leftrightarrow
\eta\in \HH_{\mu}$.
Assume that $\P[\HH_{\mu}]=1$. Then any joint measure
is almost Gibbs. This was pointed out and discussed
in the papers \cite{ku,ku1} and is apparent from
the above representation of the conditional expectation.

Let us remark that,
whenever $K$ is a translation-invariant probability
measure on the product space and $K^{\bar \si}$ is any
joint measure with marginal $K_d^{\bar \si}(d\eta)=\P(d\eta)$
we have that
$K_d(d\eta)\neq \P(d\eta) \Rightarrow h(K|K^{\bar \si})>0$.
This is clear from the monotonicity
of the relative entropy w.r.t. to the filtration.
(\cite{G}, Proposition 15.5 c). So, $h(K|K^{\bar \si})=0$
would imply that $h(K_d|\P)=0$
which again
would imply $K_d=\P$, by the classical variational principle
applied to the product measure $\P$.
So, given a joint measure $K^{\bar \si}$, the
class of interesting measures is reduced to the
ones having the same $\eta$-marginal.

\bp Suppose that $\P$ is a product measure
and that $\g^\mu$ is the above specification
for a translation-invariant joint measure $K^\mu$.
Suppose that  $\P(\HH_{\mu})=1$.
Take $K$ a translation-invariant measure  with marginal
$K_d=\P$.

Then condition C'1 holds for the measure $K$,
for the above choice of $\l$.
\ep

\bpr We have to check that
$\l(d\si^1 d\eta^1) K(d\si^2 d\eta^2) $-a.s
a configuration $\si^1_{<0} \eta^1_{<0}\si^2_{\geq 0} \eta^2_{\geq 0}$
is in $\O_{\g^{\mu}}$, where  for a configuration $\si$ we have written $\si_{<0}=(\si_x)_{x < 0}$, etc.
This is equivalent to  $ \eta^1_{<0}\eta^2_{\geq 0} \in \HH_{\mu}$ for $\P \otimes \P$-a.e. $\eta^1, \eta^2$,  since
both $\l$ and $K$ have marginal $\P$, and the later is immediate because it is a product
measure.
\epr

To illustrate the general theory of Section 3.4 we
note the following

\bc Suppose that $\P$ is a product measure
and that $\g^\mu$ is the above specification
for a translation-invariant joint measure $K^\mu$.
Suppose that  $\P(\HH_{\mu})=1$.
Take $K\in \GG_{inv}(\g^\mu)$  with marginal
$K_d=\P$.

Then condition C'2 of Theorem \ref{chp:thm6} is true and hence
\[
h(K|K^{\mu})=
\lim_{\L}\frac{1}{|\L|}\int \P(d\eta)
\log \frac{K(\si^0_{\L}\eta_{\L})}{K^{\mu} (\si^0_{\L}\eta_{\L})}=0
\]
for any translation invariant spin-configuration $\si^0$.
\ec

\subsection{Random field
Ising model: failure of the second part of the variational principle}
\bigskip

Let us now specialize to the random field
Ising model.

For all what follows we will denote
by $K^+(d\si d\eta)=\P(d\eta)\mu^+[\eta](d\si)$
the `plus-joint measure'. Here we clearly
mean by $\mu^+[\eta](d\si)=\lim_{\L\uparrow\Z^d}\mu^+[\eta](d\si)$
the random infinite-volume Gibbs measure on the Ising
spins. The limit exists for any arbitrary fixed $\eta$,
by monotonicity. Similarly we write
$K^-(d\si d\eta)=\P(d\eta)\mu^-[\eta](d\si)$.
In this situation we have
\medskip

\bp
Assume that the quenched random field Ising model
has a phase transition in the sense that
$\mu^+[\eta](\si_x=+)>\mu^-[\eta](\si_x=+)$ for $\P$-a.e. $\eta$ and for some $x \in \Z^d$.
Then the joint measures $K^+$ and $K^-$, obtained
with the same defining potential are not
compatible with the same specification.
\ep

\br
We already know by Theorem \ref{chp:thm6.2} that the relative entropy $h(K^+|K^-)$ is zero, and thus we prove here that the second part of the variational principle  is not valid in case of phase transition for the quenched random field Ising model.
\er

\br
In the so-called ``grand ensemble approach'' to disordered systems
proposed in the theoretical physics literature (going back to Morita (1964))
it is implicitly assumed
that the potential for the joint measure always exists and
does not depend on the choice of the joint measure for the
same defining potential.
Here we give a full proof that
non-unicity of the joint conditional expectations (and necessarily
of the corresponding joint potential) really does happen, 
despite of the fact that the joint
measures are always weakly Gibbs. 
It is thus an important example of a pathological
behavior in the Morita approach in a well-known
disordered system, in a translation-invariant situation. 
For a discussion of the problems of the
Morita approach with the theoretical physics community, 
see \cite{EKM,KM1,KM2}.
\er

\bpr The proof relies on the
explicit representation of proposition
\ref{kulsprop} for the conditional
expectations of $K^+$ (resp. $K^-$) in terms
of $\mu^+$ (resp. $\mu^-$).
We will show that
$\int K^+(d\xi_{x^c}) K^{-}_x(\cdot|\xi_{x^c}) \neq  K^+(\cdot)$.
Let us evaluate both sides on the event
$B:=\{\eta_x=+, \sum_{y:|y-x|=1}\si_y=0 \}$.

Using proposition \ref{kulsprop}, it is simple to see that
we have in particular for the local event $\eta_x=+$ for any
configuration $\si$ with $\sum_{y:|y-x|=1}\si_y=0$
the formula
\[
 K^+(\eta_x=+|\si_{x^c}\eta_{x^c})
=\Bigl(
1+\int \mu^+[\eta_x=-,\eta_{x^c}](d\tilde \si_x)e^{2h \tilde \si_x}
\Bigr)^{-1}=: r^{+}(\eta_{x^c})
\]
So we get that
\[
 K^+(B)=\int\P(d \tilde\eta)
\mu^+[\tilde \eta]\Bigl( \sum_{y:|y-x|=1}\si_y=0 \Bigr)
\times r^{+}(\tilde \eta_{x^c})
\]
Define $r^{-}(\eta_{x^c})$ as above, but with the
Gibbs measure $\mu^-$. Then we have
\[
\int K^+(d\xi_{x^c}) K^{-}_x(\cdot|\xi_{x^c})(B)
=\int\P(d \tilde\eta)
\mu^+[\tilde \eta]
\Bigl( \sum_{y:|y-x|=1}\si_y=0 \Bigr)
\times r^{-}(\tilde \eta_{x^c}).
\]
Now it follows from our assumption that,
for $\P$-a.e. configuration $\tilde \eta$ we have
the strict inequality
$r^{+}(\tilde \eta_{x^c})<r^{-}(\tilde \eta_{x^c})$.
But this shows that both measures give different
expectations of $B$ and finishes the claim.
\epr
\bigskip

In the following we show from
the weakly Gibbsian point of view that $K^+$ and $K^-$ have
a ``good" (rapidly decaying) almost surely
convergent translation invariant potential. This strengthens
the results in \cite{ku1}, where the a.s. absolutely convergent
potential is not translation invariant.

\begin{theorem}Assume that $d\geq 3$,
$\b$ is large enough, the random fields $\eta_x$
are i.i.d. with symmetric distribution that
is concentrated on finitely many values,
and that $h\P\eta^2_x$ is sufficiently small.

There exists an absolutely convergent potential
that is {\em translation invariant} for the plus joint measure
$K^+(d\si d\eta)$ for sufficiently low temperature and small disorder.
It decays like a stretched exponential.
\end{theorem}
\bpr
Applying the remark given after (5.5) that rely on Theorem 2.4.
of  \cite{ku1}  we have the following. \\

{\bf Fact proved in \cite{ku1}.} 

Assume that
$K^\mu(d\x)=\P(d\eta)\mu[\eta](d\si)$
is a joint measure for the random field Ising model.
Denote the disorder average of the quenched spin-spin
correlation by
\[
 c(m):=\sup_{{x,y:|x-y|=m}}
\int\P(d \eta)
\Bigl|\mu[\eta](\si_x \si_y)
-\mu[\eta](\si_x)\,
\mu[\eta](\si_y)
\Bigr|.
\]
Suppose we give ourselves any nonnegative translation invariant
function $w(A)$ giving weight
to a subset $A\sb \Z^d$.

Then there is a potential
$\bar U^{\mu}(\eta)$ on the disorder space
satisfying the decay property
\[
\sum_{A:A\ni x_0} w(A)
\int\P(d\eta)
\left |
\bar U^\mu_{A}(\eta)\right|
\leq \bar C_1+\bar C_2\sum_{m=2}^{\infty} m^{2d-1} \bar w(m) c(m)
\]
if the r.h.s. is finite.
Here $\bar w(m):= w\left(\{z\in \Z^d;z\geq  0, |z|\leq m \}\right)$
where $\geq $ denotes the lexicographic order.
$\bar C_1, \bar C_2$ are constants, depending on $\b,h$.
If $K^\mu$ is translation invariant, then
$\bar U^{\mu}(\eta)$ is translation
invariant, too.
The total potential
$U^{\rm{triv}}(\si,\eta)+\bar U^{\mu}(\eta)$
is a potential for $K^\mu$. Here
$U^{\rm{triv}}$ is a potential for the
formal Hamiltonian $-\b\sum_{<i,j>}\si_i\si_j-h\sum_{i}\eta_i\si_i
-\sum_i\log \P_0(\eta_i)$. \\

It was already stated in \cite{ku1} that we expect
a superpolynomial decay of the
quantity $c(m)$ with $m$,  when $m$ tends
to infinity. We remark first that it was already stated and proved in Bricmont and Kupiainen (1988) that
$\left|\mu[\eta](\si_x \si_y)
-\mu[\eta](\si_x)
\mu[\eta](\si_y)
\right|\leq C(\eta)e^{-C \b d(x,y)}
$
with a random constant $C(\eta)$ that is finite
for $\P$-a.e. $\eta$. The problem
is that integrability of the constant
is not to be expected.
 Unfortunately,
Bricmont and Kupiainen  (1988) do not control explicitly in their paper
the decay of the disorder average $c(m)$.
Now we will reenter their renormalization
group proof and
sketch how stretched exponential decay
is obtained for $c(m)$. Obviously,
we cannot repeat the details
of the RG analysis here. For a pedagogical exposition 
of the RG for disordered models, see 
also \cite{BoK} where the example of an interface model was treated. 

\bc{(From \cite{BK})}
There is an exponent $\a>0$ such that, for all
$m$ sufficiently large we have that
\be \label{6.54}
c(m)\leq e^{-m^\a}.
\ee
\ec

{\em Sketch of proof based on RG:}\\

For the first part we follow Bricmont and Kupiainen (1988), page 750, 8.3
`Exponential Decay of Correlations'.
Fix $x$ and $y$. We will be interested
in sending their distance to infinity.
Let us denote by $H\sb \Z^d$ the half space $H:=\{
z\in \Z^d, e\cdot z\leq a\}$ for $a>0$, where $e$ is
a fixed unit vector.
Let us denote by $\mu_H[\eta]:=\lim_{\L\uparrow H}\mu^+_{\L}[\eta]$.
By monotonicity we have for any configuration
of random fields $\eta$ that the quenched expectation
of the spin at the origin in the measure $\mu^+_a[\eta]$ is bigger
than that in the measure $\mu^+[\eta]$.

Repeating the FKG-arguments
given in the first steps of Bricmont and Kupiainen (1988), Chapter 8.3., it
is sufficient to show stretched exponential
decay
of the quantity
\[ \int\P(d\eta)\Bigl(\mu^+_H[\eta](\si_0)
-\mu^+[\eta](\si_0)\Bigr)
\]
as a function of $d(H^c,0)$ to prove (\ref{6.54}). 
As in \cite{BK} we denote by $E_H$ the ``good''
event in spin-space in all
of $\Z^d$ that there is no Peierls contour around $0$ in that touches the
complement of $H$. Then, in the same configuration $\eta$,
we have that the r.h.s. is bounded by
\[
 \mu^+_H[\eta](\si_0)
-\mu^+[\eta](\si_0)
\leq \mu^+[\eta](E_H^c).
\]
Now, we can always estimate this expectation
as a sum over probabilities of Peierls contours
\[ \mu^+[\eta](E_H^c)\leq \sum_{{\g:\inn \g\ni 0}\atop{
\inn \g \cap H^c\neq \emp}}
\mu^+[\eta](\g).
\]
The problem is that there is no uniform Peierls estimate
for all configurations of the disorder.
There is however a ``good event'' in disorder space
$G=G_H$ such that there really is a Peierls estimate
for all the ``long'' contours appearing in the above sum.
The $\P$-probability of the complement of this event is small
and controlled
(in a very-nontrivial way) by the renormalization group
construction. For $\eta\in G_H$ we really have
that
\[ \sum_{{\g:\inn \g\ni 0}\atop{
\inn \g \cap H^c\neq \emp}}
\mu^+[\eta](\g)\leq e^{-C \b d(H^c,0)}.
\]
This is stated as (8.34) in \cite{BK}.
So we have that
\[
\int \P(d\eta)\mu^+[\eta](E_H^c)\leq \P(G^c) + e^{-C \b d(H^c,0)}
\]
From the construction of the renormalization group
in Bricmont-Kupiainen we can see that $G$ is expressable
in the so-called bad fields
$\bold N^k_x(\eta)$
in the form $G=\{\eta, \bold N^k_x(\eta)=0 \,\,
\forall |x|<L,\,\, \forall k>\frac{\log d(x, H^c)}{\log L}
\}$. $L$ is a fixed finite length scale (the block-length suitably
chosen in the construction of the RG).
It appears here just as a constant.
$x\in \Z^d$ runs over sites in the lattice
and $k$ is a natural number denoting
the $k$-the application of the renormalization group
transformation.
The renormalization group gives the probabilistic control
of the form
\[
\P(\bold N^k_x(\eta)\neq 0)
\leq e^{- L^{r_1 k}}
\]
with some $r_1>0$
(this follows from \cite{BK} Lemma 1 and Lemma 2, page 563)
and so we have
\[ \P(G^c_H)\leq L^d
\sum_{k>\frac{\log d(0, H^c)}{\log L}}e^{- L^{r_1 k}}
\leq L^d e^{- d(0, H^c)^{r_2} }
\]
for $d(0, H^c)$ sufficiently large with $r_1>r_2>0$.
This proves the claim.
\epr

{\bf Acknowledgments:} 

The authors thank Aernout van Enter for interesting discussions and comments.

Christof K\"ulske,\\
WIAS, Mohrenstrasse 39,\\
 10117 Berlin, Germany.\\
  E-Mail: kuelske@wias-berlin.de\\

Arnaud Le Ny,\\
Eurandom, L.G. 1.48,\\
 TU Eindhoven, Postbus 513,\\
  5600 MB Eindhoven, The Netherlands.\\
   E-Mail: leny@eurandom.tue.nl\\

Frank Redig,\\
Faculteit Wiskunde En Informatica,\\
 TU Eindhoven, Postbus 513,\\
  5600 MB Eindhoven, The Netherlands.\\
   E-Mail:f.h.j.redig@tue.nl


\addcontentsline{toc}{section}{\bf References}
\begin{thebibliography}{14}


\bibitem{BK} J.~Bricmont, A.~Kupiainen (1988). \newblock Phase transition  in the $3d$ random field Ising model. {\em Comm. Math. Phys.} {\bf 142}:539--572.

\bibitem{BKL} J.~Bricmont, A.~Kupiainen, R.~Lefevere (1998). \newblock
  Renormalization group pathologies and the definition of Gibbs
  states. {\em Comm. Math. Phys.} {\bf 194}:359--388.

\bibitem{BoK} A.~Bovier, C.~K\"ulske (1994). \newblock A rigorous renormalization group method 
for interfaces in random media. {\em Rev. Math. Phys.}
{\bf 6}, no. 3:413--496.

\bibitem{DS} R.L. Dobrushin and S.B. Shlosman (1997). \newblock Gibbsian
  description of 'non Gibbsian' field. \newblock \emph{Russian Math.
    Surveys} {\bf 52}:285-297. Also 'Non Gibbsian' states and
  their description. \emph{Comm. Math. Phys.} {\bf 200}:125--179,
  1999.

\bibitem{EFHR}
A.C.D. van Enter, R. Fern\'andez, F. den Hollander, F. Redig (2002).
Possible loss and recovery of Gibbsianness during the stochastic
evolution of Gibbs measures. {\em Comm. Math. Phys.}
{\bf 226}:101--130.

\bibitem{EFS} A.C.D. van Enter, R. Fern{\'a}ndez and A.D. Sokal (1993).
  \newblock Regularity properties and pathologies of position-space
  renormalization-group transformations: Scope and limitations of
  {G}ibbsian theory. \newblock {\em J. Statist. Phys.} {\bf 72}:879--1167.


\bibitem{EKM} A.C.D. van Enter, C. Maes and C. K\"ulske (2000). 
\newblock Comment on: \cite{KM1}, {\em Phys. Rev. Lett.} {\bf 84}:6134.



\bibitem{EMSS} A.C.D. van Enter, C. Maes, R.H. Schonmann and S.
Shlosman (2000a). \newblock The Griffiths singularity random field.
\newblock {\em On Dobrushin's way: from probability theory to
statistical physics}, Amer. Math. Soc. Transl. Serie 2,
{\bf 198}:51--58, {\em Amer. Math. Soc. Providence, P.I}.


\bibitem{FLNR} R.~Fern\'andez, A. Le Ny and F. Redig (2002). \newblock
Variational principle and almost quasilocality for renormalized
measures.
\newblock Eurandom report 2002-032. To appear in \emph{J. Statist. Phys.} Available at http://euridice.tue.nl/\,$\tilde{}$aleny.

\bibitem{FP} R. Fern\'andez and C.-E. Pfister (1997).\newblock Global
  specifications and non-quasilocality of projections of {G}ibbs
  measures.\newblock \emph{Ann. Probab.} {\bf 25}:1284--1315.

\bibitem{G} H.O.~Georgii (1988).  \newblock Gibbs Measures and Phase
  Transitions. \newblock Walter de Gruyter (de Gruyter Studies in
  Mathematics, Vol.\ 9), Berlin-New York.


\bibitem{GP} R.B. Griffiths and P.A. Pearce (1979). \newblock Mathematical
  properties of position-space renormalization-group transformations.
  \emph{J. Statist. Phys.} {\bf 20}: 499--545.

\bibitem{I2} R.B. Israel (1986). \newblock Convexity in the theory of
lattice gases. Princeton university Press.

\bibitem{K} O.K. Kozlov (1974). \newblock Gibbs description of a system of
  random variables. \emph{Probl. Inf. Transm.} {\bf
    10}:258--265.




\bibitem{KM1} R. K\"uhn and G. Mazzeo (1994). \newblock
Critical behavior of the randomly spin diluted
                   2D Ising model: A grand ensemble approach. {\em Phys. Rev. Lett.} {\bf 73}, 2268--2271. 

\bibitem{KM2} R. K\"uhn and G. Mazzeo (2000). Reply to \cite{EKM} \newblock {\em Phys. Rev. Lett.} {\bf 84}, 6135.

\bibitem{ku}
C. K\"ulske (1999). \newblock (Non-)Gibbsianness and phase transitions in
random lattice spin models. \newblock {\em Markov Process. Related Fields} {\bf 5}, no 4:357-383.

\bibitem{ku1}
C. K\"ulske (2001). \newblock Weakly Gibbsian representation for joint
measures of quenched lattice spin models. \newblock  {\em Probab. Theory Related Fields} {\bf 119}:1--30.


\bibitem{lef} R. Lefevere (1999).\newblock Variational principle for some
  renormalized measures. \newblock {\em J. Statist. Phys.} {\bf
    95}:785--803.

\bibitem{lef2}R. Lefevere (1999a). Almost and weak Gibbsianness: A long-range
  pair-interaction example. {\em J. Statist. Phys.} {\bf 96}:109--113.

\bibitem{MRSV} C. Maes, F. Redig, S. Shlosman and A. Van Moffaert (2000).
\newblock Percolation, path large deviations and weak Gibbsianity.
\newblock {\em Comm. Math. Phys.} {\bf 209}: 517--45.

\bibitem{MRTV} C. Maes, F. Redig,  F. Takens,  A. Van Moffaert and E. Verbistky (2000).  Intermittency and weak Gibbs states,  {\em Nonlinearity} {\bf 13}, no.5:1681--1698. 


\bibitem{MRV}
C.~Maes, F.~Redig and A.~Van Moffaert (1999).
\newblock Almost Gibbsian versus Weakly Gibbsian.
\newblock \emph{Stochastic Process. Appl.} {\bf 79}:1--15.

\bibitem{MRV2}
C.~Maes, F.~Redig and A.~Van Moffaert (1999a).
\newblock The restriction of the Ising model to a Layer.
\newblock {\em J. Statist. Phys.} {\bf 96}:69--107.



\bibitem{MVDV} C. Maes and K. vande Velde (1997). \newblock Relative energies for non-Gibbsian states. \newblock
{\em Comm. Math. Phys.} {\bf 189}:277-286.


\bibitem{Mo}
T. Morita (1964). \newblock Statistical Mechanics of quenched solid solutions with application
to magnetically dilute alloys.
\newblock  {\em J. Math. Phys.} {\bf 5}:1402--1405.





\bibitem{P}
C.-E. Pfister (2002).\newblock Thermodynamical aspects of classical
lattice systems. \newblock Proceedings of the 4th brazilian school
of probability, Mambucaba, RJ, 2000. In {\em In and out of
equilibrium. Probability with a physical flavour}. Progress in
probability, Vladas Sidoravicius editor, Birkha\"user, pp 393-472.



\bibitem{Pirlot} M. Pirlot (1980). \newblock A strong variational principle for continuous spin systems. {\em J. Appl. Probab.} {\bf 17}, no. 1:47--58.



\bibitem{sepp} T. Sepp\"al\"ainen (1995). \newblock Entropy, limit theorems, and variational principles for disordered lattice systems, {\em Comm. Math. Phys.}{\bf 171}:233-277.

\bibitem{Su} W.G. Sullivan (1973). \newblock Potentials for almost Markovian
  random fields. {\em Comm. Math. Phys.} {\bf 33}:61--74.
  
\bibitem{XI} S. Xu (1997). \newblock An ergodic process of zero divergence distance from
the class of all stationary processes. {\em J. Theoret. Probab.} {\bf 11}:181-196.

\end{thebibliography}
\end{document}